\documentclass[10pt]{article}
\usepackage{bbm}
\usepackage{amsfonts}
\usepackage{amsmath}
\allowdisplaybreaks
\usepackage{amssymb}
\usepackage{fancyhdr}
\usepackage{setspace}
\usepackage{latexsym}
\usepackage{mathrsfs}
\usepackage{amsthm}

\usepackage{cite,enumitem,graphicx}
\usepackage[colorlinks=true,urlcolor=blue,
citecolor=red,linkcolor=blue,linktocpage,pdfpagelabels, bookmarksnumbered,bookmarksopen]{hyperref}
\usepackage[english]{babel}

\numberwithin{equation}{section}
\newtheorem{theorem}{Theorem}[section]
\newtheorem{proposition}[theorem]{Proposition}
\newtheorem{lemma}[theorem]{Lemma}
\newtheorem{corollary}[theorem]{Corollary}
\newtheorem{definition}[theorem]{Definition}
\newtheorem{remark}[theorem]{Remark}
\newtheorem{example}[theorem]{Example}

\newcommand{\R}{\mathbb R}

\newcommand{\e}{\varepsilon}

\newcommand{\bt}{\begin{theorem}}
\newcommand{\et}{\end{theorem}}
\newcommand{\bl}{\begin{lemma}}
\newcommand{\el}{\end{lemma}}
\newcommand{\bd}{\begin{definition}}
\newcommand{\ed}{\end{definition}}
\newcommand{\bc}{\begin{corollary}}
\newcommand{\ec}{\end{corollary}}
\newcommand{\bp}{\begin{proof}}
\newcommand{\ep}{\end{proof}}
\newcommand{\bx}{\begin{example}}
\newcommand{\ex}{\end{example}}
\newcommand{\bi}{\begin{exercise}}
\newcommand{\ei}{\end{exercise}}
\newcommand{\bo}{\begin{proposition}}
\newcommand{\eo}{\end{proposition}}
\newcommand{\br}{\begin{remark}}
\newcommand{\er}{\end{remark}}
\newcommand{\be}{\begin{equation}}
\newcommand{\ee}{\end{equation}}
\newcommand{\ba}{\begin{align}}
\newcommand{\ea}{\end{align}}
\newcommand{\bn}{\begin{enumerate}}
\newcommand{\en}{\end{enumerate}}
\newcommand{\bg}{\begin{align*}}
\newcommand{\bcs}{\begin{cases}}
\newcommand{\ecs}{\end{cases}}

\newcommand{\bean}{\begin{eqnarray*}}
\newcommand{\eean}{\end{eqnarray*}}

\def \e{\varepsilon}

\makeatletter 
\renewcommand\theequation{\thesection.\arabic{equation}}
\@addtoreset{equation}{section}
\makeatother
\numberwithin{equation}{section}

\begin{document}

\begin{center}
\textbf{Stable solution and extremal solution for fractional $p$-Laplacian}\\
\end{center}

\begin{center}
Weimin Zhang$^{1, 2}$\\

\smallskip
$^1$ School of Mathematical Sciences, Zhejiang Normal University, Jinhua 321004, P.R. China
\smallskip

$^2$ School of Mathematical Sciences,  Key Laboratory of Mathematics and Engineering Applications (Ministry of Education) \& Shanghai Key Laboratory of PMMP,  East China Normal University, Shanghai, 200241, P.R. China
\end{center}
\begin{center}
\renewcommand{\theequation}{\arabic{section}.\arabic{equation}}
\numberwithin{equation}{section}
\footnote[0]{\hspace*{-7.4mm}
{E-mail addresses: zhangweimin2021@gmail.com (W. Zhang).}}
\end{center}

\begin{abstract}
To our knowledge, this paper is the first attempt to consider the existence issue for fractional $p$-Laplacian equation: $(-\Delta)_p^s u= \lambda f(u),\; u> 0 ~\text{in}~\Omega;\;
u=0\;\text{in}~ \mathbb{R}^N\setminus\Omega$,
where $p>1$, $s\in (0,1)$, $\lambda>0$ and $\Omega$ is a bounded domain with $C^{1, 1}$ boundary. We first propose a notion of stable solution, then we prove that when $f$ is of class $C^1$, nondecreasing and satisfying $f(0)>0$ and $\underset{t\to \infty}{\lim}\frac{f(t)}{t^{p-1}}=\infty$, there exists an extremal parameter $\lambda^*\in (0, \infty)$ such that a bounded minimal solution $u_\lambda \in W_0^{s,p}(\Omega)$ exists if $\lambda\in (0, \lambda^*)$, and no bounded solution exists if $\lambda>\lambda^*$. Moreover, no $W_0^{s,p}(\Omega)$ solution exists for $\lambda > \lambda^*$ if in addition $f(t)^{\frac{1}{p-1}}$ is convex.

To handle our problems, we show a Kato-type inequality for fractional $p$-Laplacian. We show also $L^r$ estimates for the equation $(-\Delta)_p^su=g$ with $g\in W_0^{s, p}(\Omega)^*\cap L^q(\Omega)$ for $q \geq 1$, especially for $q \le \frac{N}{sp}$. We believe that these general results have their own interests. Finally, using the stability of minimal solutions $u_\lambda$, under the polynomial growth or convexity assumption on $f$, we show that the extremal function $u_* =\lim_{\lambda\to\lambda^*}u_\lambda  \in W_0^{s,p}(\Omega)$ in all dimensions, and $u^*\in L^{\infty}(\Omega)$ in some low dimensional cases. \end{abstract}

\section{Introduction}\label{s1}
We are interested in the existence problem for the following fractional $p$-Laplacian equation
\[
\leqno (P_{\lambda})~~~~~~~~~~~~~~~~~~~~~~~~~~~~~
\begin{cases}
\begin{aligned}
(-\Delta)_p^s u&= \lambda f(u) &&\text{in}~\Omega,\\
u&> 0~~~~~&&\text{in}~\Omega,\\
u&=0~~~~~~&&\text{in}~ \mathbb{R}^N\setminus\Omega,
\end{aligned}
\end{cases}
\]
where $p>1$, $s\in (0,1)$,  $\lambda>0$ and $\Omega$ is a bounded domain with $C^{1, 1}$ boundary. Here $f: \mathbb{R}\to\mathbb{R}$ satisfies
\begin{itemize}
\item[$(F_1)$]  $f$ is of class $C^1$, nondecreasing, $f(0)>0$ and $\underset{t\to \infty}{\lim}\frac{f(t)}{t^{p-1}}=\infty$.
\end{itemize}
The classical model examples are $f(u)=e^u$ and $f(u)=(1+u)^m$ with $m>p-1$. Let us begin with some well-known results with the Laplacian, the fractional Laplacian and the $p$-Laplacian.

\subsection{Laplacian case}
When $p=2$, $s=1$, $(P_\lambda)$ is the classical semilinear equation
\begin{equation}\label{2406121915}
\begin{cases}
\begin{aligned}
-\Delta u&=\lambda f(u) &&\text{in}~\Omega,\\
u&>0 &&\text{in}~ \Omega,\\
u&=0 &&\text{on}~ \partial\Omega.
\end{aligned}
\end{cases}
\end{equation}
For $f$ satisfying $(F_1)$, a well-known result is: there exists an extremal parameter $\lambda^*>0$ such that if $0<\lambda<\lambda^*$, \eqref{2406121915} admits a minimal classical solution $u_\lambda$, which is stable, increasing with respect to $\lambda$; while for $\lambda>\lambda^*$, \eqref{2406121915} has no classical solutions; see the monograph \cite{Dupaigne2011} and the survey \cite{Cabre2017}. 

When $f$ is in addition convex, Brezis {\it et al.}\cite{BC1996} introduced the notion of $L^1$-weak solution and verified that $\lambda^*$ is also a dichotomy parameter for the existence of $L^1$-weak solution. Moreover they proved that $u^*(x) = {\lim}_{\lambda\to\lambda^*} u_\lambda(x)$ is always a $L^1$-weak solution, called the {\it extremal solution}. More precisely, $u$ is said a $L^1$-weak solution to \eqref{2406121915} if 
\[
u\in L^1(\Omega), \quad f(u)\mbox{dist}(x, \partial\Omega)\in L^1(\Omega),
\]
and
\[
-\int_{\Omega}u\Delta \varphi dx=\lambda\int_{\Omega}f(u) \varphi dx\quad\mbox{$\forall\;\varphi\in C^2_{c}(\Omega)$.}
\]
A solution $u$ to \eqref{2406121915} is called to be stable if 
\[
\int_{\Omega}|\nabla \varphi |^2 dx\ge \lambda\int_{\Omega}f'(u) \varphi^2 dx\quad\mbox{$\forall\;\varphi\in H^1_0(\Omega)$.}
\]

The regularity of the extremal solution has attracted many attentions. When $f(u)=e^u$, we know that $u^*$ is bounded if $N\le 9$ in \cite{CR1975}, while $u^*=-2\log{{|x|}}$ if $N\ge 10$ and $\Omega=B_1$ in \cite{JL1972}.  Crandall and Rabinowitz \cite{CR1975}, Mignot and Puel \cite{MP1980} studied the case $f(u)=(1+u)^m$ with $m>1$. They showed that $u^*$ is bounded  if
\begin{align*}
N< G_1(m, 2):=2+\frac{4m}{m-1}+4\sqrt{\frac{m}{m-1}}.
\end{align*}
Joseph and Lundgren \cite{JL1972} used phase plane techniques to show that $u^*$ is unbounded when $\Omega$ is a ball and $N \geq G_1(m, 2)$. We can remark that $G_1(m, 2)>10$ for all $m>1$.

Later, Brezis and V\'azquez \cite{BV1997, Brezis2003} raised some open questions about the regularity of $u^*$. In particular, for all convex nonlinearities $f$ satisfying $(F_1)$ and smooth bounded domains $\Omega$, they asked whether $u^*$ belongs always to $H^1_0(\Omega)$, and whether $u^*$ is always bounded for $N\le 9$.

Twenty years having passed, many partial answers were established, see for instance \cite{Nedev2000, CC2006, Cabre2010, Cabre2019, Villefas2013, CR2013, YZ2002}. In 2020, Cabr\'e, Figalli, Ros-Oton and Serra \cite{CFRS2020} completely answered these questions on the regularity of $u^*$. They showed that if $\Omega$ is a $C^3$ domain, then $u^*\in H_0^1(\Omega)$ in all dimensions and $u^*\in L^{\infty}(\Omega)$ for $N\le 9$.
\subsection{$p$-Laplacian case}
For $p>1$, $s=1$, $(P_\lambda)$ becomes the quasilinear $p$-Laplacian equation
\begin{equation}\label{2306182235}
\begin{cases}
\begin{aligned}
-\Delta_p u&=\lambda f(u) &&\text{in}~\Omega,\\
u&>0 &&\text{in}~\Omega,\\
u&=0 &&\text{on}~ \partial\Omega.
\end{aligned}
\end{cases}
\end{equation}
Cabr\'e and Sanch\'on \cite{CS2007} provided a notion of stability to \eqref{2306182235} for general nondecreasing $f$, a solution $u\in W_0^{1, p}(\Omega)$ of \eqref{2306182235} is said to be stable if
\[
\int_{\{\nabla u(x)\neq 0\}}\Big[ |\nabla u|^{p-2}|\nabla \varphi|^2 + (p-2)|\nabla u|^{p-4}(\nabla u \nabla \varphi)^2\Big] dx\ge \lambda\int_{\Omega}f'(u)\varphi^2dx \quad\mbox{$\forall\;\varphi\in \mathcal{A}_u(\Omega)$},
\]
where the space of test functions $\mathcal{A}_u(\Omega)$ is chosen as
\[
\mathcal{A}_u(\Omega):=
\begin{cases}
W_{0}^{1,p}(\Omega), &\mbox{ if } p\ge 2,\\
\big\{\varphi\in W_0^{1, p}(\Omega):|\varphi |\le C|u|,\;  |\nabla\varphi |\le C|\nabla u| \mbox{ for some } C>0\big\},& \mbox{ if } 1<p<2.
\end{cases}
\]
With the assumption $(F_1)$, they derived a dichotomy result for the existence of bounded solutions to \eqref{2306182235}. The study on the regularity of $u^*$ with $p$-Laplacian began with Garc\'ia Azorero, Peral Alonso
and Puel \cite{GP1992, GP1994}, they proved that for $f(u)=e^u$, $u^*\in L^{\infty}(\Omega)$ if
\begin{equation*}
N< p+\frac{4p}{p-1},
\end{equation*}
 while $u^*=-p\log{{|x|}}$, $\lambda^*=p^{p-1}(N-p)$ if $N\ge p+\frac{4p}{p-1}$ and $\Omega=B_1$. Cabr\'e and Sanch\'on \cite{CS2007} considered power-like nonlinearities $f\in C^1(\R_+)$ satisfying $(F_1)$ and 
\begin{equation*}\label{2311092040}
0\le f(t)\le C(1+t)^m \;\mbox{ in $\R_+$ for some $C>0$, $m > p-1$;} \quad \mbox{and }\; \underset{t\to \infty}{\lim\inf}\frac{f'(t)t}{f(t)}\ge m.
\end{equation*}
They obtained then $u^*\in W_0^{1, p}(\Omega)$ for all $N\ge 1$ and $u^*\in L^{\infty}(\Omega)$ for
\begin{equation}\label{2311041433}
N<G_1(m, p):=\frac{p}{p-1}\left(1+\frac{mp}{m-p+1}+2\sqrt{\frac{m}{m-p+1}}\right).
\end{equation}
Here \eqref{2311041433} is optimal, since when $\Omega=B_1$ and $f(u)=(1+u)^m$, in \cite{CS2007, Ferrero2004} it was shown that for all $N\ge G_1(m, p)$, $m>p-1$,
$$u^*=|x|^{\frac{-p}{m-(p-1)}}-1 \quad \mbox{and}\quad \lambda^*=\left(\frac{p}{m- p+1}\right)^{p-1}\left[N-\frac{mp}{m-p+1}\right].$$

If moreover
\begin{itemize}
\item[$(F_2)$]~~~~~~~~~~~~~~$
\begin{cases}
\begin{aligned}
&(f-f(0))^{\frac{1}{p-1}}\mbox{ is convex} &&\text{for}~p\ge 2,\\
&f\mbox{ is convex} &&\text{for}~1<p<2,\\
\end{aligned}
\end{cases}
$
\end{itemize}
the $W_0^{1,p}(\Omega)$ regularity for $u^*$ was obtained by Sanch\'on \cite[Theorem 1]{Sanchon_PA2007} if $p > 2$ and $N<p+\frac{p^2}{p-1}$; by Castorina and Sanch\'on \cite[Theorem 1.5]{CS2015} in any dimension if $p>2$ and $\Omega$ is a strictly convex domain. When $p\in (1, 2]$, Sanch\'on \cite[Theorem 2]{Sanchon07} proved $u^*\in W_0^{1,p}(\Omega)$ in all dimensions.

The boundedness of $u^*$ is also studied under the assumptions $(F_1)$-$(F_2)$ in low dimensions. Sanch\'on \cite{Sanchon07, Sanchon_PA2007} proved that for general bounded smooth domains $\Omega$, $u^*$ is bounded whenever
\begin{equation}\label{2311091916}
N<p+\frac{p}{p-1} \;\;\text{if}~p\ge 2; \quad  N\le p+\frac{2p}{p-1}\left(1+\sqrt{2-p}\right) \;\;\text{if}~1<p<2.
\end{equation}
Cabr\'e, Miraglio and Sanch\'on \cite{CMS2022} studied also the case of strictly convex domains $\Omega$.
For more details, we refer to \cite{CMS2022} and references therein.

\subsection{Fractional Laplacian case}
As far as we know, only a few results have been established about extremal solution in fractional Laplacian problem. Ros-Oton and Serra \cite[Proposition 1.2]{RS2014} proved a dichotomy result of bounded solutions about the fractional semilinear equation ($0 < s < 1$)
\begin{equation}\label{2313032137}
\begin{cases}
\begin{aligned}
(-\Delta)^s u&=\lambda f(u) &&\text{in}~\Omega,\\
u&>0 &&\text{in}~\Omega,\\
u&=0 &&\text{in}~\mathbb{R}^N\backslash\Omega.
\end{aligned}
\end{cases}
\end{equation}
They showed that if $f$ is convex, then $u^*\in L^\infty(\Omega)$ whenever $N<4s$. Following the approach in \cite{Sanchon07}, they also established the boundedness of $u^*$ in dimensions $N<10s$ assuming that $f\in C^2$ and the limit $\underset{t\to\infty}{\lim}\frac{f(t)f''(t)}{f'(t)^2}$ exists. When $\Omega$ is convex, the Pohozaev identity is used in \cite{RS2014} to deduce $u^*\in H^s(\mathbb{R}^N)$ for all $N\ge 1$ and $s\in(0, 1)$.

For the exponential case $f(u)=e^u$, when $\Omega$ is symmetric and convex with respect to all the coordinate directions, Ros-Oton \cite{Ros-Oton2014} obtained that $u^*\in L^{\infty}(\Omega)$ whenever
\[
\frac{\Gamma^2(\frac{N+2s}{4})}{\Gamma^2(\frac{N-2s}{4})}<\frac{\Gamma(\frac{N}{2})\Gamma(1+s)}{\Gamma^2(\frac{N-2s}{2})}.
\]
Here $\Gamma$ is Euler's Gamma function. In particular, the extremal solution $u^*$ is bounded for all $s\in (0,1)$ whenever $N\le 7$.

More generally, when $\Omega$ is a ball and $f$ is a $C^2$ function satisfying $(F_1)$, Sanz-Perela \cite{San-Perela2018} proved $u^*\in L^{\infty}(\Omega)$ for  $2\le N< 2s+4+2\sqrt{2(s+1)}$.
For $s=\frac12$,  Cabr\'e and Sanz-Perela \cite{CS2022} showed that if moreover $f$ and $\Omega$ are smooth and convex,  then $u^*\in L^{\infty}(\Omega)$ for $1\le N\le 4$.

\subsection{Our main results}
As mentioned, we will consider the existence problem to the fractional $p$-Laplacian problem $(P_\lambda)$, we aim to establish some dichotomy type results and also study the regularity of the extremal function $u^*$. The crucial point is to establish some uniform estimates for minimal solutions $u_\lambda$, where the stability is a key ingredient. As explained, the stability stands for the non-negativeness of the second variation of associated energy functional.  To our best knowledge, it seems that there is no work in this direction for the fractional $p$-Laplacian case. 

Motivated by the mentioned works, we are interested in providing a notion of stable solution for the following more general problem
\begin{equation}\label{2306182230}
(-\Delta)_p^s u =g(x, u) \;\;\text{in}~\Omega,\quad
u =0 \;\;\text{in}~ \mathbb{R}^N\setminus\Omega
\end{equation}
where $p>1$, $s\in (0,1)$. We introduce first the weak formulation of $(-\Delta)_p^s$. Denote the Gagliardo seminorm by
\[
[u]_{s,p}:=\left( \int_{\mathbb{R}^{2N}}{\frac{|u(x)-u(y)|^p}{|x-y|^{N+sp}}}dxdy\right)^{1/p}.
\]
Let
\[
W^{s,p}(\mathbb{R}^N):=\{u\in L^p(\mathbb{R}^N): [u]_{s,p}<\infty\}
\]
be endowed with the norm
\[
\|u\|_{W^{s,p}(\mathbb{R}^N)}:=(\|u\|_p^p+[u]_{s,p}^p)^{1/p}.
\]
Here and after $\|\cdot\|_q$ denotes the usual norm of $L^q(\mathbb{R}^N)$ for $q \geq 1$. Let
\[
W_0^{s,p}(\Omega):=\left\{u\in W^{s,p}(\mathbb{R}^N):  u=0 \; \text{a.e.~in}~\mathbb{R}^N \setminus \Omega\right\},
\]
equivalently renormed with $\|u\|_{W_0^{s,p}(\Omega)}=[u]_{s,p}$ (see \cite[Theorem 7.1]{Di12}), it is well-known that $W_0^{s,p}(\Omega)$ is a uniformly convex Banach space for any $p \in (1, \infty)$. When $\Omega$ has a continuous boundary, $C_c^{\infty}(\Omega)$ is dense in $W_0^{s,p}(\Omega)$, see \cite[Theorem 6]{FSV2015}. $(-\Delta)_p^s$ can be variationally regarded as an operator from $W_0^{s,p}(\Omega)$ into its dual space $W_0^{s, p}(\Omega)^*$ as follows,
\begin{align}
\label{FpL}
\langle (-\Delta)_p^s u, v \rangle=\int_{\mathbb{R}^{2N}}{\frac{J_u(x,y)(v(x)-v(y))}{|x-y|^{N+sp}}}dxdy, \quad \forall\; v\in W_0^{s,p}(\Omega),
\end{align}
where 
\begin{align}
\label{Jw}
J_w(x,y)=|w(x)-w(y)|^{p-2}(w(x)-w(y)), \quad \forall\; w\in W_0^{s,p}(\Omega).
\end{align}
In the last decade, many studies for fractional $p$-Laplacian were realized, see for instance, the comparison principle \cite[Proposition 2.10]{IMS2016} (see also \cite[Lemma 9]{LL2014}), the maximum principle \cite[Theorem A.1]{BF2014}, $L^{\infty}$-estimate \cite{BP2016}, Hopf's lemma \cite{DQ2017} and H\"older regularity theory \cite{IMS2016}.
\begin{definition}
\label{Def1.1}
We say $u\in W_0^{s,p}(\Omega)$ is a weak solution (respectively subsolution, supersolution) to \eqref{2306182230} if $g(x, u)\in W_0^{s, p}(\Omega)^*$ and
\begin{equation*}
\langle (-\Delta)_p^s u, v \rangle=(\text{respectively}~\le,\, \ge) \int_{\Omega} g(x, u)vdx
\end{equation*}
holds for all $v \in W_0^{s,p}(\Omega)$ with $v\ge 0$. By $g(x, u)\in W_0^{s, p}(\Omega)^*$, we mean that 
for any $v \in W_0^{s,p}(\Omega)$, there holds ${g}(x, u)v \in L^1(\Omega)$ and
\[
\Big|\int_{\Omega} {g}(x, u)vdx \Big|\le C\|v\|_{W_0^{s,p}(\Omega)}.
\]
\end{definition}

\begin{remark}
 Without confusion, sometimes we omit the word {\it weak} for simplicity; some other times, to emphasize the $W_0^{s,p}(\Omega)$ regularity, we also call {\it weak solution} as {\it $W_0^{s,p}(\Omega)$ solution}.
\end{remark}

Formally speaking, the energy functional associated to \eqref{2306182230} is
\begin{equation}\label{2307061350}
I(u)=\frac1{p} \|u\|_{W_0^{s, p}(\Omega)}^p-\int_{\Omega}\left(\int_0^{u}g(x, t) dt\right) dx,\quad \forall\; u\in W_0^{s,p}(\Omega).
\end{equation}
Notice that when $1<p<2$ or the reaction term $g$ possesses Sobolev supercritical growth, the functional $I$ might be not twice differentiable, or even not well-defined in $W_0^{s, p}(\Omega)$. Therefore, when $1<p<2$, it's more subtle to give a comprehensive definition of stability. Inspired by \cite{CS2007}, we choose the test function space at $u\in W_0^{s,p}(\Omega)$ as follows.
\[
\begin{aligned}
\mathcal{T}_{u}(\Omega):=
\begin{cases}
W_0^{s,p}(\Omega)&\text{if}~p\ge 2;\\
\{\varphi\in W_0^{s,p}(\Omega): \;|\varphi(x)-\varphi(y)|\le C_\varphi|u(x)-u(y)|~\mbox{a.e.~for some }C_\varphi>0\}&\text{if}~1<p< 2.
\end{cases}
\end{aligned}
\]
Note that the above constant $C_\varphi$ could depend on $u$ and $\varphi$. If $1<p<2$, for any $\varphi\in \mathcal{T}_u(\Omega)$, since $u(y)=\varphi(y)=0$ for $y\in \mathbb{R}^N\backslash\Omega$, we actually have $|\varphi(x)|\le C_\varphi|u(x)|$. 

\begin{definition}\label{2406152227}
Assume that $p>1$, $s\in (0,1)$, and $g(x, t)$ is differentiable and nondecreasing in $t$ for a.e. $x\in\Omega$. Let $u\in W_0^{s,p}(\Omega)$ be a weak solution to \eqref{2306182230}. We say that $u$ is stable if
\begin{equation}\label{2306201531}
(p-1)\int_{\mathbb{R}^{2N}}{\frac{|u(x)-u(y)|^{p-2}(\varphi(x)-\varphi(y))^2}{|x-y|^{N+ps}}}dxdy-\int_{\Omega}{g_t}(x, u)\varphi^2 dx\ge 0,
\end{equation}
for all $\varphi \in\mathcal{T}_{u}(\Omega)$.
\end{definition}
\begin{remark}
We select $\mathcal{T}_{u}(\Omega)$ as space of test functions to make the nonlocal integral in \eqref{2306201531} well-defined. Note that the second integral in \eqref{2306201531} is well defined with values in $[0, \infty]$ since $g_t(x, u)\ge 0$. For $1 < p < 2$, comparing to the choice of Cabr\'e and Sanch\'on \cite{CS2007}, we don't require the control of $\nabla\varphi$.
\end{remark}

Our first main result for $(P_\lambda)$ is a dichotomy claim for existence of bounded solution.
\begin{theorem}\label{2306211344}
Assume that $p>1$, $s\in (0,1)$ and $f$ satisfies $(F_1)$. Then there exists $\lambda^*\in(0, \infty)$ such that $(P_\lambda)$ has a bounded minimal solution $u_\lambda\in W_0^{s, p}(\Omega)$ for $\lambda\in (0, \lambda^*)$, which is also a stable solution. Moreover, $u_\lambda$ is nondecreasing in $\lambda$, and no bounded solution exists for $\lambda>\lambda^*$. 
\end{theorem}
By virtue of the super-sub solution method, we can get a family of minimal solutions $u_\lambda$, but we need some new argument to claim the stability of $u_\lambda$. The reason is that $(-\Delta)_p^s$ is a nonlinear and nonlocal operator, Theorem \ref{2306211344} could not be obtained directly as in fractional Laplacian case or $p$-Laplacian case. In order to handle the stability, we make use of a cut-off functional. More precisely, consider a positive weak supersolution $\overline{u}$ of \eqref{2306182230}, let
\[
\widetilde{g}(x, t):=
\begin{cases}
g(x, \overline{u}(x)),\;\; t> \overline{u}(x),\\
g(x, t), \;\;~~~~ 0\le t\le \overline{u}(x),
\end{cases}
\]
and
\begin{equation}\label{2311301928}
\widetilde{I}(u)=\frac{1}{p}\|u\|_{W_0^{s,p}(\Omega)}^p-\int_{\Omega}{\widetilde{G}}(x, u) dx,\quad \forall\; u\in W_0^{s,p}(\Omega),
\end{equation}
where
\[
\widetilde{G}(x, u)=\int_{0}^{u^+(x)}\widetilde{g}(x, t) dt;\quad u^+ =\max\{u, 0\}.
\]
We will check that $\widetilde{I}$ has a global minimizer $\widehat u$ in $W_0^{s, p}(\Omega)$ (see Proposition \ref{2308061322}), which is a stable solution of
\begin{equation*}
(-\Delta)_p^s u = \widetilde{g}(x, u^+) \;\;\text{in}~\Omega,\quad
u =0 \;\;\text{in}~ \mathbb{R}^N\setminus\Omega.
\end{equation*}
Using comparison principle \cite[Proposition 2.10]{IMS2016}, we can claim finally $\widehat u$ is a solution of \eqref{2306182230}. Moreover, with the condition $g(x, 0)\not\equiv 0$, there exist sufficiently many test functions in the set $\{\varphi\in \mathcal{T}_{\widehat u}(\Omega): 0\le \varphi \le \widehat u\}$, which yields the stability of $\widehat u$.

Come back to the problem $(P_\lambda)$. To verify the stability of minimal solution $u_\lambda$ with the definition \ref{2406152227}, we set $\overline{u}=u_\lambda$ as a weak supersolution and work with the functional $\widetilde{I}$ in \eqref{2311301928}, consequently the corresponding minimizer $\widehat u$ is identical to $u_\lambda$. It is worthy to remark that our method in Section \ref{2312021503} can simplify the proof for $p$-Laplacian case \cite[Theorem 1.4(i)]{CS2007}.

\medskip
Furthermore, if in addition $f(t)^{\frac{1}{p-1}}$ is convex, $\lambda^*$ is also the dichotomy parameter for $W_0^{s,p}(\Omega)$ solution.
\begin{theorem}\label{2312030111}
Assume that $p>1$, $s\in (0,1)$ and $f$ satisfies $(F_1)$. If $f(t)^{\frac{1}{p-1}}$ is convex, then $(P_\lambda)$ has no $W_0^{s,p}(\Omega)$ solution for $\lambda>\lambda^*$ where $\lambda^*$ is the critical value given in Theorem \ref{2306211344}.
\end{theorem}

When $p\in (1, 2)$, even for $\varphi\in C_c^{\infty}(\Omega)$, one cannot claim in general
\begin{equation}\label{2312031549}
\begin{aligned}
\int_{\mathbb{R}^{N}}\frac{|\varphi(x)-\varphi(y)|^{p-2}(\varphi(x)-\varphi(y))}{|x-y|^{N+sp}}dy\in L_{loc}^{\infty}(\Omega),
\end{aligned}
\end{equation}
see \cite[Remark 2.14]{IMS2016}. Thus for $u\in W_0^{s, p}(\Omega)$, we could not get easily the Kato-type inequality
\begin{equation}\label{2312271147}
 (-\Delta)_p^s \Psi(u) \le|\Psi'(u)|^{p-2}\Psi'(u)  (-\Delta)_p^s u
\end{equation}
for general convex and Lipschitz $\Psi\in C^{2}(\mathbb{R})$ satisfying $\Psi(0)=0$. However, 
this kind of argument is important to handle Theorem \ref{2312030111}. To overcome this difficulty, we establish a Kato-type inequality slightly weaker than \eqref{2312271147}. Consider the equation
\begin{equation}\label{2306190123}
(-\Delta)^s_p u =g(x) \;\;\text{in}~\Omega;\quad
u =0\;\; \text{in}~ \mathbb{R}^N\setminus\Omega.
\end{equation}
\begin{proposition}\label{2308041528}
Assume that $p>1$, $0<s<1$, $g\in W_0^{s, p}(\Omega)^*\cap L^1(\Omega)$ and $u\in W_0^{s, p}(\Omega)$ is a weak solution of \eqref{2306190123}. Let $\Psi\in C^{2}(\mathbb{R})$ be convex and  {Lipschitz}, $\Psi(0)=0$. If in addition either $\Psi''\in L^{\infty}$ or $g$ has a constant sign, there holds
\begin{equation}\label{2304242313}
\begin{aligned}
\langle (-\Delta)_p^s \Psi(u), \psi\rangle \le \int_{\Omega}|\Psi'(u)|^{p-2}\Psi'(u) g \psi dx\quad \mbox{for all}\; \psi\in W_0^{s, p}(\Omega) \cap L^{\infty}(\Omega), \psi\geq 0.
\end{aligned}
\end{equation}
\end{proposition}
\begin{remark}
In virtue of Proposition \ref{2308041528}, we will establish in Proposition \ref{2312141015} that any weak solution to \eqref{2306190123} satisfies 
\begin{equation}\label{2312271423}
\langle (-\Delta)_p^s|u|, \psi \rangle\le \int_{\Omega}|g|\psi dx \quad \mbox{ for all } \psi\in W_0^{s,p}(\Omega)\cap L^{\infty}(\Omega) , \psi\ge 0.
\end{equation}
Here the test function space is chosen as $W_0^{s,p}(\Omega)\cap L^{\infty}(\Omega)$, since we actually cannot claim $|g|\in W_0^{s, p}(\Omega)^*$. 
\end{remark}

Next, we are concerned with the understanding of the extremal function
\begin{equation}\label{2311151118}
u^*(x) :=\underset{\lambda\to \lambda^*}{\lim}u_\lambda(x),
\end{equation}
where $\{u_\lambda\}$ is given in Theorem \ref{2306211344}. A natural question is whether $u^*$ possesses $W_0^{s, p}(\Omega)$ or $L^{\infty}(\Omega)$ regularity. Notice that even if $u^*\in W_0^{s, p}(\Omega)$, to show that $u^*$ is a solution of $(P_{\lambda^*})$, one needs further to prove $f(u^*)\in W_0^{s, p}(\Omega)^*$. 

Consider first \eqref{2306182230} with power-like nonlinearities, i.e. $g$ satisfies
\begin{equation}\label{2306202233}
0\le g(x, t)\le C(1+t)^m;\quad \text{$g$ is nondecreasing with respect to $t$}
\end{equation}
for $t\ge 0$, a.e.~$x\in \Omega$; and
\begin{equation}\label{2306211610}
\underset{t\to \infty}{\lim\inf}\frac{g_t(x,t)t}{g(x, t)}\ge m \quad \mbox{uniformly for a.e. $x\in\Omega$.}
\end{equation}
Let $m>p-1$, denote
\begin{equation}\label{2311092123}
G_s(m, p):=\frac{sp}{p-1}\left(1+\frac{mp}{m-p+1}+2\sqrt{\frac{m}{m-p+1}}\right).
\end{equation}
The following theorem gives the uniform $W_0^{s, p}(\Omega)$ or $L^\infty(\Omega)$ bound for stable solutions to \eqref{2306182230}.
\begin{theorem}\label{2306211611}
Assume that $p>1$, $s\in (0,1)$, $m>p-1$, $g$ satisfies \eqref{2306202233}-\eqref{2306211610}.  There exists $C$ depending only on $s$, $m$, $p$, $\Omega$ and $g$ such that if $u\in W_0^{s,p}(\Omega)$ is a non-negative stable solution of \eqref{2306182230}, then
\[
\|u\|_{W_0^{s,p}(\Omega)}\le C.
\]
If in addition $N<G_s(m, p)$, there holds $\|u\|_{\infty}\le C$.
\end{theorem}
Remark that when $u^*$ is bounded, $u^*\in C^{\alpha}(\overline{\Omega})$ for some $\alpha\in (0, s]$ thanks to \cite[Theorem 1.1]{IMS2016}. Since $u^*$ in \eqref{2311151118} is the pointwise limit of $u_\lambda$, its regularity can be derived easily from uniform estimates of $u_\lambda$.
\begin{theorem}\label{2307200048}
Assume that $p>1$, $s\in (0,1)$, $m>p-1$ and $f$ satisfies $(F_1)$,
\begin{align}
\label{newThm1.9}
{\liminf_{t\to\infty}\frac{f'(t)t}{f(t)} \geq m, \quad \limsup_{t\to\infty}\frac{f(t)}{t^m} < \infty.}
\end{align}
Then the following assertions hold:
\begin{itemize}
\item[\rm (i)] $u^*\in W_0^{s,p}(\Omega)$ is a stable solution of $(P_{\lambda^*})$;
\item[\rm (ii)] If moreover $N<G_s(m, p)$, then $u^*\in L^{\infty}(\Omega)$.
\end{itemize}
\end{theorem}

\begin{remark}
Clearly, $G_s(m, p)>sp+\frac{4sp}{p-1}$ for all $m>1$. Hence $u^* \in L^\infty(\Omega)$ provided $f$ satisfies $(F_1)$, \eqref{newThm1.9} and $N\le sp+\frac{4sp}{p-1}$.
\end{remark}

The proof of Theorem \ref{2306211611} is similar to \cite[Theorem 1.2]{CS2007}, however some new obstacles occur due to the nonlocal nature of $(-\Delta)_p^s$. A main argument is to understand the regularity of equation \eqref{2306190123} provided $g\in W_0^{s,p}(\Omega)^*\cap L^q(\Omega)$ for $q\ge 1$. 
The following estimates have their own interests for the study of fractional $p$-Laplacian problem.
\begin{proposition}\label{2306271234}
Let $s\in (0, 1)$, $p>1$ and $\Omega$ be a bounded domain. Assume that $g\in  W_0^{s,p}(\Omega)^*\cap L^q(\Omega)$ for some $q\ge 1$, and $u$ is the weak solution to \eqref{2306190123}. Then we have the following regularity results: There exist constants $C=C(\Omega, p, s, q, r)$ such that
{\begin{itemize}
\item[\rm (i)] If $q=1$, for any $0< r< \frac{N(p-1)}{N-sp}$ with $N > sp$; or any $r \geq 1$ with $N = sp$,
\begin{equation}\label{2309152339}
\|u\|_r\le C\|g\|_1^\frac{1}{p-1};
\end{equation}
\item[\rm (ii)] If $q\in (1, \frac{N}{sp})$ {\rm(}hence $N > sp${\rm )}, 
\begin{equation}\label{2306242254-1}
\|u\|_{\frac{N(p-1)q}{N-spq}}\le C\|g\|_q^\frac{1}{p-1};
\end{equation}
\item[\rm (iii)] If $q=\frac{N}{sp} > 1$, for any $\alpha_0<\frac{N}{N-sp}\Big(\frac{p^pS}{\|g\|_{\frac{N}{sp}}}\Big)^\frac{1}{p-1}$,
\begin{equation}\label{2312142133}
\int_{\Omega}e^{\alpha_0 |u|} dx<\infty,
\end{equation}
where $S$ is the best Sobolev embedding constant defined in \eqref{eq2406131933} below;
\item[\rm (iv)] If $q>\frac{N}{sp}\geq 1$,
\begin{equation}\label{2306242254}
\|u\|_{\infty}\le C\|g\|_q^\frac{1}{p-1};
\end{equation}
\item[\rm (v)] If $q\ge 1 > \frac{N}{sp}$, $\alpha=\frac{sp-N}{p}$,
\[
\|u\|_{C^{\alpha}(\overline{\Omega})}\le C\|g\|_q^\frac{1}{p-1}.
\]
\end{itemize}}
\end{proposition}
In this paper, when $N>sp$, $S$ denotes the best Sobolev constant defined by
\begin{equation}\label{eq2406131933}
S:=\underset{u\in D^{s,p}(\mathbb{R}^N) \backslash \{ 0\}}{\text{inf}}\frac{[u]_{s,p}^p}{|u|_{p_s^*}^{p}},
\end{equation}
where $
D^{s,p}(\mathbb{R}^N):=\big\{u\in L^{p^*_s}(\mathbb{R}^N): [u]_{s,p}< \infty\big\}$ and $p_s^*=\frac{Np}{N-sp}$.
\begin{remark}
Proposition \ref{2306271234} extends \cite[Proposition 1.4]{RS2014} and \cite{LPPS2015} with fractional Laplacian and \cite[Lemma 3.1]{CS2007} with $p$-Laplacian. Our proof of Proposition \ref{2306271234} is elementary by choosing suitable test functions in \eqref{2312271423}, which seems to be simpler than the approach of \cite{CS2007}. Furthermore,  \cite{BP2016} and \cite{IMS2016} respectively derived $L^{\infty}$-estimate for $q>\frac{N}{sp}$ and H\"older estimate for $q=\infty$. Notice also that the solution is continuous when $q>\frac{N}{sp}$, seeing \cite[Theorem 3.13]{BP2016} and \cite[Theorem 1.5]{KMS2015},  and the H\"older continuity was proven in \cite{BLS2018} if in addition $p\ge 2$. 
\end{remark}

The following results consider general nonlinearities $f\in C^2(\R_+)$ satisfying $(F_1)$-$(F_2)$, including for example the exponential growth case as $f(u)=e^u$. 
\begin{theorem}\label{2309111457}
Let $p\ge 2$ and $f$ satisfies $(F_1)$, $(F_2)$. If $N<sp(1+\frac{p}{p-1})$, then $u^*$ is a {stable} $W_0^{s,p}(\Omega)$ solution of $(P_{\lambda^*})$. Moreover, let $p' = \frac{p}{p-1}$ be the conjugate of $p$, then
\[
\begin{aligned}
u^*\in
\begin{cases}
L^{\infty}(\Omega), &\mbox{ if }N<spp',\\
L^{q}(\Omega) \;\;\mbox{ for all }1\le q<\infty, &\mbox{ if }N=spp',\\
L^{q}(\Omega) \;\;\mbox{ for all }1\le q<\frac{N(p-1)^2}{N(p-1)-sp^2}, &\mbox{ if }N>spp'.\\
\end{cases}
\end{aligned}
\]
\end{theorem}
\begin{theorem}\label{2307201716}
Assume that $p>1$, $s\in (0,1)$ and $f\in C^2(\R_+)$ satisfies $(F_1)$, $(F_2)$. If 
\begin{equation}\label{2308071234}
\tau:=\underset{t\to\infty}{\lim}\frac{f(t)f''(t)}{f'(t)^2} \;\mbox{exists and }\; \tau > \frac{p-2}{p-1},
\end{equation} 
then $u^*$ is a {stable}  $W_0^{s,p}(\Omega)$ solution of $(P_{\lambda^*})$. Moreover,  $u^* \in L^{\infty}(\Omega)$ when $N<sp+\frac{4sp}{p-1}$.
\end{theorem}

\begin{remark}
It's not difficult to see that the extremal solution $u^*$ in Theorems \ref{2307200048},  \ref{2309111457} and \ref{2307201716} is the minimal solution of $(P_{\lambda^*})$.
\end{remark}

Theorem \ref{2307200048} extends \cite[Theorem 1.4(ii)]{CS2007} to fractional $p$-Laplacian case. As mentioned in \eqref{2311041433}, for the boundedness of $u^*$, $G_1(m, p)$ is optimal in the $p$-Laplacian case with $f(u) = (1+u)^m$. However, the approach in \cite{CS2007, Ferrero2004} cannot work with the fractional $p$-Laplacian, even for $f(u)=(1+u)^m$ on the unit ball.  Theorems \ref{2309111457}, \ref{2307201716} extend \cite{Sanchon07, Nedev2000, Sanchon_PA2007} to the fractional $p$-Laplacian. Again the bound $sp+\frac{4sp}{p-1}$ is consistent to the optimal result for the $p$-Laplacian equation with $f(u)=e^u$. 

However, regardless of the values of $s$ and $m, p > 1$, $u^*$ is showed in \cite{San-Perela2018} bounded whenever $2\le N\le 6$ with $f(u) = (1+u)^m$, \cite{Ros-Oton2014} obtained that $u^*\in L^{\infty}(\Omega)$ for $N\le 7$  with $f(u) = e^u$. So $G_s(m, p)$ or $sp+\frac{4sp}{p-1}$ is not always optimal for the fractional $p$-Laplacian problem $(P_\lambda)$, the general answer seems unknown even for $p=2$.


\medskip
Our paper is organized as follows. In Sections \ref{2408111106}-\ref{2312021506}, we respectively establish the Kato-type inequality (Proposition \ref{2308041528}) and the regularity result (Proposition \ref{2306271234}) for fractional $p$-Laplacian. Section \ref{2312021503} is dedicated to the existence of stable solution where we show dichotomy results for bounded solutions and $W_0^{s,p}(\Omega)$ solutions. In Section \ref{2312021519}-\ref{2312021527}, we consider the regularity of extremal function with power-like nonlinearities and convex nonlinearities, and prove Theorem \ref{2307200048} and Theorems \ref{2309111457}, \ref{2307201716} respectively.

\section{Kato inequality for fractional $p$-Laplacian}\label{2408111106}
Here we prove the Kato type inequality for fractional $p$-Laplacian, i.e. Proposition \ref{2308041528}. 
The following inequality can refer to \cite[Lemmas A.1]{BP2016}.
\begin{lemma}\label{2312192211}
Let $1<p<\infty$ and $h:\mathbb{R}\to \mathbb{R}$ be a $C^1$ convex function. For $\ell\ge 0$, we set
\[
\xi_{p, \ell}(t)=\left(\ell+|t|^2\right)^{\frac{p-2}{2}}t,\quad t\in\mathbb{R}
\]
then
\[
\begin{aligned}
|a-b|^{p-2}(a-b)&\left[\alpha\xi_{p, \ell}(h'(a))-\beta \xi_{p, \ell}(h'(b))\right]\\
\ge& \left(\ell(a-b)^2+(h(a)-h(b))^2\right)^{\frac{p-2}{2}}(h(a)-h(b))(\alpha - \beta)
\end{aligned}
\]
for any $a \ne b\in\mathbb{R}$ and any $\alpha,\beta \ge 0$.
\end{lemma}

\begin{proof}[\bf Proof of Proposition  \ref{2308041528}]
Let first $g\ge 0$ and $g_k(x)=\min\{k,  g(x)\}$. By \cite[Theorem 3.1]{BP2016}, there exists a sequence $\{u_k\}\subset L^{\infty}(\Omega)\cap W_0^{s, p}(\Omega)$ such that $(-\Delta)_p^s u_k=g_k$. From comparison principle \cite[Proposition 2.10]{IMS2016}, we see that
\[
\| u_k\|_{W_0^{s, p}(\Omega)}^p =\int_{\Omega}g_k u_k dx\le \int_{\Omega} gu dx =\|u\|_{W_0^{s, p}(\Omega)}^p.
\]
Hence $\{u_k\}$ is bounded in $W_0^{s, p}(\Omega)$. Up to a subsequence, we can assume that $u_k\to u$ weakly in $W_0^{s,p}(\Omega)$ and pointwisely a.e.~in $\Omega$.

For $p\ge 2$, take any $0 \le \psi \in W^{s, p}_0(\Omega)\cap L^\infty$, one can check that  $\varphi=\psi |\Psi'(u_k)|^{p-2}\Psi'(u_k)\in W_0^{s, p}(\Omega)$. Using $\varphi$ as test function to \eqref{2306190123} and apply Lemma \ref{2312192211} with
\[
h=\Psi,\;\;\ell=0,\;\; a=u_k(x),\;\; b=u_k(y),\;\; \alpha=\psi(x)\; \mbox{ and  }\;  \beta =\psi(y),
\]
there holds 
\begin{equation}\label{2312140959}
\begin{aligned}
\int_{\mathbb{R}^{2N}}\frac{J_{\Psi(u_k)}(x, y)(\psi(x)-\psi(y))}{|x-y|^{N+sp}}dxdy\le \int_{\Omega}g_k|\Psi'(u_k)|^{p-2}\Psi'(u_k)\psi dx,
\end{aligned}
\end{equation}
where $J$ is defined in \eqref{Jw}. Owing to the fact that
\[
\Theta_k(x, y):=\frac{J_{\Psi(u_k)}{(x, y)}}{|x-y|^{\frac{(N+sp)(p-1)}{p}}}\mbox{ is bounded in } L^{\frac{p}{p-1}}(\mathbb{R}^{2N})
\]
and pointwisely converges to
\[
\Theta(x, y):= \frac{J_{\Psi(u)}{(x, y)}}{|x-y|^{\frac{(N+sp)(p-1)}{p}}} \mbox{ for a.e.~in } \; \mathbb{R}^{2N}.
\]
Up to a subsequence, we may assume that $\Theta_k$ weakly converges to $\Theta$ in $L^{\frac{p}{p-1}}(\mathbb{R}^{2N})$,
hence
\[
\int_{\mathbb{R}^{2N}}\frac{J_{\Psi(u_k)}(x, y)(\psi(x)-\psi(y))}{|x-y|^{N+sp}}dxdy\longrightarrow \int_{\mathbb{R}^{2N}}\frac{J_{\Psi(u)}(x, y)(\psi(x)-\psi(y))}{|x-y|^{N+sp}}dxdy.
\]
Furthermore, since
\[
|\Psi'(u_k)|^{p-1}|g_k\psi|\le \|\Psi'\|_{\infty}^{p-1}|g|\psi\in L^1(\R^N),
\]
by Lebesgue's theorem, passing to the limit in \eqref{2312140959} as $k\to\infty$, we get \eqref{2304242313}.

For $1<p<2$, we can reset the test function as $\varphi=\psi \xi_{p, \varepsilon}(\Psi'(u_k))$ with $\varepsilon > 0$,
and use Lemma \ref{2312192211} with the same parameters as above but $\ell=\varepsilon$. It follows that
\[
\begin{aligned}
\int_{\mathbb{R}^{2N}}&\frac{\big[\varepsilon(u_k(x)-u_k(y))^2+(\Psi(u_k(x))-\Psi(u_k(y)))^2\big]^{\frac{p-2}{2}}\big[\Psi(u_k(x))-\Psi(u_k(y))\big] \big[\psi(x)-\psi(y)\big]}{|x-y|^{N+sp}}dxdy\\
& \le \int_{\Omega}g_k|\Psi'(u_k)|^{p-2}\Psi'(u_k)\psi dx.
\end{aligned}
\]
By weak convergence argument, letting $k\to\infty$ and then $\varepsilon\to 0$, we get again \eqref{2304242313}.

Moreover, if we have no sign information on $g$, but $\|\Psi''\|_\infty < \infty$, we can choose directly the test function by $\varphi=\psi |\Psi'(u)|^{p-2}\Psi'(u)$ for $p\ge 2$ or $\varphi=\psi (\varepsilon+|\Psi'(u)|^2)^{\frac{p-2}{2}}\Psi'(u)$ for $p\in (1, 2)$. Thus we can carry out the above steps to claim the same conclusion.
\end{proof}

The following version of Kato inequality is well-known for classical cases, but seems to be new for fractional $p$-Laplacian.
\begin{proposition}\label{2312141015}
Let $g \in W_0^{s, p}(\Omega)^*\cap L^1(\Omega)$ and $u$ be a weak solution to \eqref{2306190123}. Then 
\begin{equation}\label{2306241807}
\int_{\mathbb{R}^{2N}}\frac{J_{|u|}(x,y)(\psi(x)-\psi(y))}{|x-y|^{N+sp}}dxdy\le \int_{\Omega}|g|\psi dx
\end{equation}
for all $\psi\in W_0^{s, p}(\Omega)\cap L^{\infty}(\Omega)$ with $\psi\ge 0$.
\end{proposition}
\begin{proof}[\bf Proof]
 Set $f_{\varepsilon}(t)=(\varepsilon^2+t^2)^{\frac12}-\varepsilon$, 
one can check that $f_{\varepsilon}$ satisfies all assumptions of Proposition  \ref{2308041528}.
It results that for any non-negative $\psi\in W_0^{s, p}(\Omega)\cap L^{\infty}(\Omega)$, there holds
\begin{equation}\label{2312140956}
\int_{\mathbb{R}^{2N}}\frac{J_{f_{\varepsilon}(u)}(x, y)(\psi(x)-\psi(y))}{|x-y|^{N+sp}}dxdy\le \int_{\Omega}g|f'_{\varepsilon}(u)|^{p-2}f'_{\varepsilon}(u)\psi dx.
\end{equation}
Clearly, $f_{\varepsilon}$ is 1-Lipschitz, hence $|f'_{\varepsilon}(u)|^{p-1}|g\psi|\le |g|\psi$ and $|f_{\varepsilon}(u(x))-f_{\varepsilon}(u(y))|\le |u(x)-u(y)|.$ By Lebesgue's theorem and passing to the limit as $\varepsilon\to 0$ in \eqref{2312140956}, we get \eqref{2306241807}. \end{proof}

\section{Integral estimates for fractional $p$-Laplacian}\label{2312021506}
This section is devoted to $L^r$ estimate for the fractional $p$-Laplacian equation \eqref{2306190123}. Since $(-\Delta)_p^s$ is a homeomorphism from $W_0^{s,p}(\Omega)$ to $W_0^{s,p}(\Omega)^*$, \eqref{2306190123} has a solution $u\in W_0^{s,p}(\Omega)$ if and only if $g\in W_0^{s,p}(\Omega)^*$. We will assume $g\in L^q(\Omega)\cap W_0^{s,p}(\Omega)^*$ and consider a priori estimate for solution to \eqref{2306190123}. Here we need another elementary inequality in \cite{BP2016}.
\begin{lemma}\label{2312192311}
Let $1<p<\infty$ and $h:\mathbb{R}\to \mathbb{R}$ be an increasing function. Define
\[
H(t)=\int_0^th'(\sigma)^\frac1p d\sigma,\quad t\in\mathbb{R},
\]
there holds
\[
|a-b|^{p-2}(a-b)\left(h(a)-h(b)\right)\ge |H(a)-H(b)|^p\quad \mbox{for all $a, b\in \mathbb{R}$.}
\]
\end{lemma}

\begin{proof}[\bf Proof of Proposition \ref{2306271234}] Look at the five cases (i)-(v) respectively.

(i) {\sl Case $q=1,\, N\ge sp$.}  Consider first $N>sp$. For $k>0$, setting $T_k(t):=\max\{-k, \min\{k, t\}\}$, we test \eqref{2306190123} by $T_k(u)$, 
\[
\|T_k(u)\|_{W_0^{s,p}(\Omega)}^p\le \langle (-\Delta)_p^s u, T_k(u)\rangle = \int_{\Omega}gT_k(u) dx.
\]
By the Sobolev inequality, we have
\[
S\|T_k(u)\|_{p_s^*}^p\le k\int_{\Omega}|g|dx,
\]
which implies
\begin{equation*}
kS^{\frac{1}{p-1}}|\{x\in \Omega: |u(x)|>k\}|^{\frac{N-sp}{N(p-1)}}\le \|g\|_1^{\frac1{p-1}}.
\end{equation*}
For any $0< r< \frac{N(p-1)}{N-sp}$ and $\ell > 0$, we have
\begin{equation}\label{2309182002}
\begin{aligned}
\|u\|_r^r& =r\int_0^{\infty}t^{r-1}|\{x\in \Omega: |u(x)|>t\}|dt\\
& \le r\int_0^\ell t^{r-1}|\{x\in \Omega: |u(x)|>t\}|dt+r\int_\ell^{\infty}t^{r-1}|\{x\in \Omega: |u(x)|>t\}|dt\\
& \le |\Omega|\ell^r +r\int_\ell^{\infty}t^{r-1}|\{x\in \Omega: |u(x)|>t\}|dt\\
& \le |\Omega|\ell^r +r\int_\ell^{\infty}t^{r-\frac{N(p-1)}{N-sp}-1}dt\times \underset{k\ge 0}{\sup} \Big(k^{\frac{N(p-1)}{N-sp}}{|\{x\in \Omega: |u(x)|>k\}|}\Big)\\
& \le |\Omega|\ell^r + \frac{r}{\frac{N(p-1)}{N-sp}-r}S^{\frac{-N}{N-sp}}\|g\|_1^{\frac{N}{N-sp}}\ell^{r-\frac{N(p-1)}{N-sp}}.
\end{aligned}
\end{equation}
Then \eqref{2309182002} yields \eqref{2309152339} by choosing
\[
\ell^{\frac{N(p-1)}{N-sp}} = \frac{r\|g\|_1^{\frac{N}{N-sp}}}{S^{\frac{N}{N-sp}}|\Omega|}.
\]
When $N=sp$, we replace the above Sobolev inequality by any embedding $W_0^{s,p}(\Omega)\subset L^r(\Omega)$ with $r > 1$.

\smallskip
(ii) {\sl Case $1<q<\frac{N}{sp}$.} The proof mainly follows the idea of \cite[Theorem 3.1]{BP2016}. For $M, \beta, \delta >0$, choose the test function as
\[
\psi=(u_M+\delta)^{\beta}-\delta^\beta, \quad \mbox{where }\; u_M=\min\{|u|, M\}.
\]
Clearly, $\psi$ belongs to $W_0^{s,p}(\Omega)\cap L^{\infty}(\Omega)$. The Kato inequality in Proposition \ref{2312141015} implies that
\begin{equation}\label{2309112155}
\int_{\mathbb{R}^{2N}}\frac{J_{|u|}(x,y)\big[(u_M+\delta)^{\beta}(x)-(u_M+\delta)^{\beta}(y)\big]}{|x-y|^{N+sp}}dxdy\le \int_{\Omega}|g|\big[(u_M+\delta)^{\beta} - \delta^\beta\big]dx.
\end{equation}
Applying Lemma \ref{2312192311} with $t \ge 0$, $h(t)=(\min\{t, M\}+\delta)^\beta$,
\[
H(t)=\int_0^th'(\sigma)^{\frac1p} d\sigma=\frac{\beta^\frac{1}{p} p}{\beta+p-1}\left[(\min\{t, M\}+\delta)^{\frac{\beta+p-1}{p}} - \delta^{\frac{\beta+p-1}{p}}\right],
\]
and using \eqref{2309112155}, we have
\begin{equation}\label{2306241956}
\frac{\beta p^p}{(\beta+p-1)^p}\big\|(u_M+\delta)^{\frac{\beta+p-1}{p}}-\delta^{\frac{\beta+p-1}{p}}\big\|^p_{W_0^{s,p}(\Omega)}\le \int_{\Omega}|g|\big[(u_M+\delta)^{\beta} - \delta^\beta\big] dx.
\end{equation}
By Sobolev and H\"older's inequality, \eqref{2306241956} yields that
\begin{equation}\label{2306242016}
\begin{aligned}
\frac{\beta p^p}{(\beta+p-1)^p}S\|(u_M+\delta)^{\frac{\beta+p-1}{p}}-\delta^{\frac{\beta+p-1}{p}}\|^p_{p_s^*}&\le \int_{\Omega}|g|\big[(u_M+\delta)^{\beta} - \delta^\beta\big] dx \le \|g\|_q\|(u_M+\delta)^\beta\|_{q'},
\end{aligned}
\end{equation}
where $q' = \frac{q}{q-1}$. The left hand side of \eqref{2306242016} can be lower bounded as follows
\begin{equation}\label{2306242026}
\|(u_M+\delta)^{\frac{\beta+p-1}{p}}-\delta^{\frac{\beta+p-1}{p}}\|^p_{p_s^*} \ge \frac{1}{2^{p-1}}\|(u_M+\delta)^{\frac{\beta+p-1}{p}}\|_{p_s^*}^p-\delta^{\beta+p-1}|\Omega|^{\frac{N-sp}{N}}.
\end{equation}
Combining \eqref{2306242016} with \eqref{2306242026}, $\frac{(\beta+p-1)^p}{\beta p^p}\ge 1$ (see \cite[Lemma A.5]{BP2016}) and $u_M \ge 0$, we obtain
\begin{align}
\label{2306242157}
\begin{split}
\|(u_M+\delta)^{\frac{\beta+p-1}{p}}\|_{p_s^*}^p & \le 2^{p-1}\delta^{\beta+p-1}|\Omega|^{\frac{N-sp}{N}}+\frac{2^{p-1}(\beta+p-1)^p}{\beta p^pS}\|g\|_q\|(u_M+\delta)^\beta\|_{q'}\\
& \le \frac{(2\delta)^{p-1}(\beta+p-1)^p}{\beta p^p}\left(|\Omega|^{1-\frac1{q'}-\frac{sp}{N}}+\frac{\|g\|_q}{S\delta^{p-1}}\right)\|u_M+\delta\|^\beta_{\beta q'}.
\end{split}
\end{align}
As $1<q<\frac{N}{sp}$, take now
\[
\beta=\frac{N(p-1)(q-1)}{N-spq} \;\; \mbox{so that} \;\; \frac{\beta+p-1}{p}p_s^*=\beta q'.
\]
Then \eqref{2306242157} gives
\begin{equation}\label{2306242212}
\|u_M+\delta\|_{\beta q'}\le \frac{2(\beta+p-1)^\frac{p}{p-1}}{(\beta p^p)^{\frac{1}{p-1}}}\left(\delta^{p-1}|\Omega|^{1-\frac1{q'}-\frac{sp}{N}}+\frac{\|g\|_q}{S}\right)^{\frac{1}{p-1}}.
\end{equation}
Passing $\delta\to 0$ and $M \to \infty$ in \eqref{2306242212}, the estimate \eqref{2306242254-1} holds true.

(iii) {\sl Case $q=\frac{N}{sp}>1$.} Here we use the test function $\psi=e^{\alpha u_M}-1$ with $u_M = \min\{|u|, M\}$ and $\alpha > 0$. By Proposition \ref{2312141015}, we have
\begin{equation}\label{2312162218}
\int_{\mathbb{R}^{2N}}\frac{J_{|u|}(x,y)(e^{\alpha u_M}(x)-e^{\alpha u_M}(y))}{|x-y|^{N+sp}}dxdy\le \int_{\Omega}|g|(e^{\alpha u_M}-1)dx.
\end{equation}
Applying Lemma \ref{2312192311} to $h(t)=e^{\alpha \min\{t, M\}}$,
\[
H(t)=\int_0^th'(\sigma)^{\frac1p} d\sigma=\frac{p}{\alpha^{1-\frac1p}}\left(e^{\frac{\alpha}p \min\{t, M\}}-1\right),
\]
and using \eqref{2312162218}, we get
\begin{equation}\label{2312170009}
\frac{p^p}{\alpha^{p-1}}\|e^{\frac{\alpha}p u_M}-1\|^p_{W_0^{s,p}(\Omega)}\le \int_{\Omega}|g|(e^{\alpha u_M}-1)dx \le \|g\|_{\frac{N}{sp}}\|e^{\frac{\alpha}{p} u_M}\|_{p_s^*}^p.
\end{equation}
By the Sobolev inequality, 
\begin{equation}\label{2312170115}
\frac{p^p}{\alpha^{p-1}}S\|e^{\frac{\alpha}p u_M}-1\|^p_{p_s^*}\le \|g\|_{\frac{N}{sp}}\|e^{\frac{\alpha}{p} u_M}\|_{p_s^*}^p.
\end{equation}
For any $\varepsilon\in (0, 1)$, the left hand side in \eqref{2312170115} can be estimated by 
\begin{equation}\label{2312170116}
\begin{aligned}
\frac{p^p}{\alpha^{p-1}}S\|e^{\frac{\alpha}p u_M}-1\|^p_{p_s^*}& \ge\frac{p^p}{\alpha^{p-1}}S\left(\|e^{\frac{\alpha}p u_M}\|_{p_s^*}-|\Omega|^{\frac1{p_s^*}}\right)^{p}\ge \frac{p^p}{\alpha^{p-1}}S(1-\varepsilon)\|e^{\frac{\alpha}p u_M}\|^p_{p_s^*}-C_\varepsilon,
\end{aligned}
\end{equation}
where $C_\varepsilon>0$ depends only on $\varepsilon, p, s, \Omega$.
Let $\alpha > 0$ satisfy $\alpha^{p-1} < \frac{Sp^p}{ \|g\|_{\frac{N}{sp}}}$, we fix $\varepsilon\in (0, 1)$ such that
\begin{equation}\label{2408091416}
\alpha^{p-1} < \frac{Sp^p(1-\varepsilon)}{ \|g\|_{\frac{N}{sp}}}.
\end{equation}
Combining \eqref{2312170115}-\eqref{2408091416}, there holds then
\begin{equation}\label{2312170125}
\left(\frac{p^p}{\alpha^{p-1}}S(1-\varepsilon)-\|g\|_{\frac{N}{sp}}\right)\|e^{\alpha u_M}\|_{\frac{N}{N-sp}}\le C_{\varepsilon}.
\end{equation}
Passing $M\to\infty$ in \eqref{2312170125}, we get \eqref{2312142133} with $\alpha_0=\frac{N\alpha}{N-sp}$.

\smallskip
(iv) {\sl Case $q>\frac{N}{sp}\ge 1$.} The subcase $N>sp$ is proven by \cite[Theorem 3.1]{BP2016}.  When $N=sp$, starting from  \eqref{2306241956}, and using Moser’s iteration technique as in \cite[Theorem 3.1]{BP2016}, we can claim \eqref{2306242254}.

\smallskip
(v) {\sl Case $q\ge 1>\frac{N}{sp}$.} The embedding $W_0^{s,p}(\Omega)\subset C^{\alpha}(\overline{\Omega})$ (see \cite[Theorem 8.2]{Di12}) implies that $L^1(\Omega)\subset W_0^{s, p}(\Omega)^*$. If $g\in L^1(\Omega)$, then $\|u\|_{C^{\alpha}(\overline{\Omega})}\le C\|g\|_1^{\frac1{p-1}}$, the proof is done.
\end{proof}

When $g$ has the polynomial growth
\begin{equation}\label{23112135}
|g(x, t)|\le C(1+t)^m,\quad m>0,
\end{equation}
using Proposition \ref{2306271234}, the bootstrap argument can lead to higher regularity result for solutions to \eqref{2306182230}. In particular, we claim
\begin{lemma}\label{2306271914}
Assume that $g: \Omega\times\mathbb{R}\to \mathbb{R}$ satisfies \eqref{23112135} for some $m>p-1$. Let $u$ be a weak solution to \eqref{2306182230}. If $g(x, u)\in W_0^{s, p}(\Omega)^*\cap L^{q_0}(\Omega)$ with
\begin{equation}\label{2405191528}
q_0 > \left(1-\frac{p-1}{m}\right)\frac{N}{sp},\quad q_0 \ge 1,
\end{equation}
then $\|u\|_{\infty}<C$
where $C>0$ depends only on $p$, $s$, $m$, $q_0$, $\Omega$ and $\|g(x, u)\|_{q_0}$.
\end{lemma}
\begin{proof}
If $q_0> \frac{N}{sp}$, Proposition \ref{2306271234} yields the boundedness of $u$.  If $q_0 = \frac{N}{sp}$, by Proposition \ref{2306271234}, $u \in L^r(\Omega)$ for any $r > 1$, so is $g(x, u)$ seeing \eqref{23112135}, we get again the boundedness of $u$. 

\smallskip
Consider now $1<q_0< \frac{N}{sp}$, it follows from Proposition \ref{2306271234} that
\begin{equation}\label{2306271250}
\|u\|_{r_1}\le C\|g(x, u)\|_{q_0}^\frac{1}{p-1}\quad \text{for}~ r_1:=\frac{N(p-1)q_0}{N-spq_0}.
\end{equation}
Seeing \eqref{23112135} and \eqref{2306271250}, we get $\|g(x, u)\|_{q_1}\le C(1+\|u\|_{r_1}^m)$ for $q_1:= \frac{r_1}{m}$. Moreover, \eqref{2405191528} means $q_1 > q_0$. By induction, whenever $q_n < \frac{N}{sp}$, there exist 
\begin{equation}\label{2307181614}
r_{n+1}=\frac{N(p-1)q_n}{N-spq_n}, \quad q_{n+1}= \frac{r_{n+1}}{m} = \frac{N(p-1)q_n}{m(N-spq_n)}
\end{equation}
and $u\in L^{r_{n+1}}(\Omega)$, $g(x, u)\in L^{q_{n+1}}(\Omega)$. In that case, it is easy to see that $q_{n+1} > q_n$. 

We claim that there exists $n_0$ such that $q_{n_0} \geq \frac{N}{sp}$. Otherwise $\lim_{n\to\infty}q_n = \ell < \infty$. Passing to the limit in \eqref{2307181614}, there holds 
\[
\ell =\frac{N(p-1) \ell}{m(N-sp\ell)}, \quad \mbox{namely }\; \ell =\left(1-\frac{p-1}{m}\right)\frac{N}{sp}.
\]
That contradicts \eqref{2405191528} since $\ell > q_0$. Thus $q_{n_0} \geq \frac{N}{sp}$ for some $n_0$, then we conclude $u \in L^\infty(\Omega)$ with the estimate as before since $g(x, u) \in L^{q_{n_0}}(\Omega)$. The case $q_0 = 1$ can be treated by the same procedure, so we omit the details. 
\end{proof}

\section{Dichotomy results for $(P_\lambda)$}\label{2312021503}
In this section, we show Theorems \ref{2306211344}-\ref{2312030111}. The main tool is the super-sub solution method, which can be stated as follows.
\begin{lemma}\label{Super_Sub_Le}
Let $g: \Omega\times\mathbb{R}\to \mathbb{R}$ be a Carath\'{e}odory function satisfying that $g(x, t)$ is continuous and nondecreasing in $t$ for a.e.~$x\in\Omega$. Let $\underline{u}\in W_0^{s,p}(\Omega)$ be a subsolution to \eqref{2306182230}, and $\overline{u}\in W_0^{s,p}(\Omega)$ be a supersolution to \eqref{2306182230} such that \begin{equation*}\label{2309180211}
\underline{u}\le \overline{u} \quad \mbox{and} \quad |g(\cdot, \underline{u}(\cdot))|,~|g(\cdot, \overline{u}(\cdot))|\in {W_0^{s,p}(\Omega)^*}.
\end{equation*}
Then there is a solution $u$ of \eqref{2306182230} such that $\underline{u}\le u \le \widetilde{u}$ for any supersolution $\widetilde{u}$ of \eqref{2306182230} satisfying $\widetilde{u}\ge \underline{u}$.
\end{lemma}
\begin{proof}[\bf Proof]
Let $u_0:=\underline{u}$, and by induction we define $u_{j+1}$ $(j\ge 0)$ as solutions of
\begin{equation}\label{2208222050}
(-\Delta)_p^s u_{j+1} =g(x, u_j) \;\;\text{in}~\Omega,\quad 
u_{j+1} =0 \;\; \text{in}~ \mathbb{R}^N\setminus\Omega.
\end{equation}
By comparison principle, $\underline{u}\le u_j \le u_{j+1} \le \overline{u}$. Consequently,
test \eqref{2208222050} by $u_{j+1}$, 
\[
\begin{aligned}
\lVert u_{j+1}\rVert_{W_0^{s,p}(\Omega)}^{p}=\int_{\Omega} g(x, u_j)u_{j+1} dx \le \int_{\mathbb{R}^N}  \max\{|\underline{u}|, |\overline{u}|\}(|g(x, \underline{u})|+|g(x, \overline{u})|) dx < \infty.
\end{aligned}
\]
Here we used Definition \ref{Def1.1}  and $\max\{|\underline{u}|, |\overline{u}|\} \in W_0^{s, p}(\Omega)$. Therefore $\{u_j\}$ is nondecreasing and bounded in $W_0^{s,p}(\Omega)$, hence $u_j$ weakly converges to some $u \in W_0^{s,p}(\Omega)$, $u_j\rightarrow u$ a.e. in $\mathbb{R}^N$. So $u$ is a weak solution to \eqref{2306182230} and $\underline{u}\le u \le \overline{u}$. As we can replace $\overline{u}$ by any supersolution $\widetilde{u} \ge \underline{u}$, there holds $u\le \widetilde{u}$.
\end{proof}

The next lemma is an elementary but useful weak convergence argument, whose proof is omitted.
\begin{lemma}\label{2310231013}
Assume that $\{u_j\}\subset W_0^{s,p}(\Omega)$ weakly converges to some $u$ in $W_0^{s,p}(\Omega)$, and $h$ is a Lipschitz continuous function with $h(0)=0$. Then $h(u_j) \mbox{ weakly converges to } h(u) \mbox{ in } W_0^{s,p}(\Omega)$.
\end{lemma}

\subsection{Dichotomy result for bounded solutions}
Let $\overline{u}$ be a supersolution of \eqref{2306182230}. Using the global minimum of $\widetilde{I}$ given in \eqref{2311301928}, we will derive a stable solution to \eqref{2306182230}. More precisely,
\begin{proposition}\label{2308061322}
Assume that $ g(x, 0)\ge 0$ is nontrivial in $\Omega$, $g(x, t)$ is nondecreasing in $t$ for a.e.~$x\in\Omega$. Let $\overline{u}\ge 0$ be a supersolution of \eqref{2306182230}. Then $\widetilde{I}$ in \eqref{2311301928} attains its global infimum in $W_0^{s, p}(\Omega)$ at some $\widehat u\le \overline{u}$, which is a weak solution of \eqref{2306182230}. If in addition $g$ is of class $C^1$, and there exists $h\in L^1(\Omega)$ such that
\begin{equation}\label{2307172130}
\sup_{0 \le w\le \widehat u}g_t(x, w)\times \widehat u^2(x)\le h(x),\quad{a.e.~in }\; \Omega
\end{equation}
then $\widehat u$ is stable.
\end{proposition}
\begin{proof}[\bf Proof]
By Definition \ref{Def1.1}, as $g(x, \overline u) \in W_0^{s,p}(\Omega)^*$, $0 \leq \widetilde{G}(x, u) \leq g(x, \overline u)u^+ \in L^1(\Omega)$, $\widetilde I$ is well defined.  

We claim that $u \mapsto \|\widetilde{G}(x, u)\|_{L^1}$ is weakly continuous in $W_0^{s,p}(\Omega)$. Indeed, let $\{v_j\}$ converge weakly to $v\in W_0^{s,p}(\Omega)$, 
\begin{equation}\label{2307171613}
\left| \widetilde{G}(x, v_j)-\widetilde{G}(x, v) \right| =\Big| \int_{v^+}^{v_j^+}\widetilde{g}(x, t)dt \Big| \le  {g}(x, \overline{u})\left|v_j^+-v^+\right|
\end{equation}
By Lemma \ref{2310231013}, $v_j^+ - v^+ \rightharpoonup 0$ hence $|v_j^+-v^+|\rightharpoonup 0$ in $W_0^{s,p}(\Omega)$. Since $g(x, \overline{u})\in W_0^{s,p}(\Omega)^*$, $\widetilde{I}$ is weakly lower semi-continuous. Moreover,
\[
\Big|\int_{\Omega} \widetilde{G}(x, v) dx \Big|\le \int_\Omega {g}(x, \overline{u})v^+ dx \le \|{g}(x, \overline{u})\|_{W_0^{s,p}(\Omega)^*}\|v^+\|_{W_0^{s,p}(\Omega)} \le \|{g}(x, \overline{u})\|_{W_0^{s,p}(\Omega)^*}\|v\|_{W_0^{s,p}(\Omega)},
\]
which readily yields that $\widetilde{I}$ is coercive since $p>1$.

Therefore, there exists a minimum point $\widehat u\in W_0^{s,p}(\Omega)$ for $\widetilde{I}$. It is easy to see that $\widehat u$ is a solution of
\begin{equation}\label{2306182256}
(-\Delta)_p^s u = \widetilde{g}(x, u^+) \;\;\text{in}~\Omega,\quad 
u =0\;\;\text{in}~ \mathbb{R}^N\setminus\Omega.
\end{equation}
By the comparison principle \cite[Proposition 2.10]{IMS2016} (see also \cite[Lemma 9]{LL2014}), there holds $\widehat u\le \overline{u}$. Using the maximum principle  \cite[Theorem A.1]{BF2014}, $\widehat u>0$ in $\Omega$. Therefore $\widetilde{g}(x, \widehat u)={g}(x, \widehat u)$, in other words $\widehat u$ is a solution of \eqref{2306182230}.

Take $g_{k_0}(x)=\min\{k_0, g(x,0)\} \in L^\infty(\Omega)$ with $k_0$ large enough such that $g_{k_0}\ne 0$. Let $u_0\in W_0^{s,p}(\Omega)$ resolve  $(-\Delta)_p^s u_0= g_{k_0}$. By the H\"older regularity result in \cite{IMS2016} and the maximum principle, we get $u_0\in C^{\alpha}(\overline{\Omega})$ for some $\alpha\in (0, s]$ and $u_0>0$ in $\Omega$. Moreover, $\widehat u \ge u_0 > 0$ by the comparison principle.

Suppose first $p\ge 2$. For any $\varphi\in C_c^{0, 1}(\Omega)$ with $\varphi\ge 0$, there exists $\gamma_0$ such that $\widehat u-\gamma_0\varphi \ge 0$. As $\widehat u$ is a global minimum point of $\widetilde{I}$, then
\begin{equation}\label{2307172148}
K(\gamma) := \widetilde{I}(\widehat u-\gamma\varphi)\ge \widetilde{I}(\widehat u).
\end{equation}
By \eqref{2307172130}, we check easily that $K$ is $C^2$ in $[0, \gamma_0)$ and 
\begin{align}\label{2307172126}
\begin{split}
K''(\gamma) &= (p-1)\int_{\mathbb{R}^{2N}}{\frac{|\widehat u(x)-\widehat u(y)-\gamma(\varphi(x)-\varphi(y))|^{p-2}(\varphi(x)-\varphi(y))^2}{|x-y|^{N+ps}}}dxdy\\
& \quad  -\int_{\Omega}g_t(x, \widehat u-\gamma\varphi)\varphi^2 dx.
\end{split}
\end{align}
We get \eqref{2306201531} with $\widehat u$ by $K''(0) \geq 0$ for all $\varphi\in C_c^{0, 1}(\Omega)$ with $\varphi\ge 0$. 
Let $\varphi\in C_c^{0, 1}(\Omega)$ (without sign constrain), as $|\varphi|\in C_c^{0, 1}(\Omega)$, \eqref{2306201531} with $|\varphi|$ implies obviously \eqref{2306201531} with $\varphi$. 

Since $C_c^\infty(\Omega)$ is dense in $W_0^{s,p}(\Omega)$ (see  \cite[Theorem 6]{FSV2015}),  \eqref{2306201531} holds in $W_0^{s, p}(\Omega)$, hence $\widehat u$ is stable.

\smallskip
Let now $1<p<2$. As above, we need only to consider $\varphi\in \mathcal{T}_{\widehat u}(\Omega)$ with $\varphi\ge 0$. By $|\varphi(x)-\varphi(y)|\le C_\varphi|\widehat u(x)-\widehat u(y)|$ for some $C_\varphi >0$, there holds $0 \le \varphi \le C_\varphi \widehat u$. So for $\gamma_0 > 0$ small enough, $\widehat u - \gamma_0\varphi\ge 0$
 and
 \begin{equation}\label{2309112015}
 \frac12 |\widehat u(x)-\widehat u(y)|\le |\widehat u(x)-\widehat u(y)- \gamma(\varphi(x)-\varphi(y))|\le 2|\widehat u(x)-\widehat u(y)|,\;\; \forall\, \gamma\in (0, \gamma_0).
 \end{equation}
Again we check that $K$ given by \eqref{2307172148} satisfies \eqref{2307172126} in $[0, \gamma_0)$. The proof is completed.
\end{proof}

Now we are in position to prove Theorem \ref{2306211344}.
\begin{proof}[\bf Proof of Theorem \ref{2306211344} completed]
It is clear that for any $\lambda>0$, $0$ is a subsolution of $(P_\lambda)$. We shall find a supersolution of $(P_\lambda)$ by classical idea. Let $\xi$ satisfy
\begin{equation*}\label{2306181420}
(-\Delta)_p^s \xi = 1 \;\; \text{in}~\Omega, \quad \xi=0 \;\; \text{in}~ \mathbb{R}^N\setminus\Omega.
\end{equation*}
By the H\"older regularity result in \cite{IMS2016}, $\xi \in C^{\alpha}(\overline{\Omega})$ for some $\alpha\in (0, s]$. We get readily $(-\Delta)_p^s \xi=1>\lambda f(\xi)$ for $\lambda > 0$ but small enough. In virtue of Lemma \ref{Super_Sub_Le}, there exists a minimal solution $0\le u_\lambda\le \xi$, hence $u_\lambda\in L^{\infty}(\Omega)$.

For $\lambda_0>0$, if $(P_{\lambda_0})$ has a bounded solution $u_{0}$, then $u_{0}$ is a supersolution to $(P_\lambda)$ for all $\lambda\in (0, \lambda_0)$. Applying Lemma \ref{Super_Sub_Le}, $(P_\lambda)$ has a bounded minimal solution $u_\lambda$ for all $\lambda\in (0, \lambda_0)$. We claim that
\[
\lambda^*:=\sup\{\lambda>0: (P_\lambda) \text{~has~a~bounded~solution}\} < \infty.
\]
In fact, let $\lambda_1$ be the first eigenvalue of $(-\Delta)^s_p$, it is known (see \cite[Theorem 4.1]{BP2016}) that $\lambda_1$ is an isolated eigenvalue, so there is $\mu >\lambda_1$, which is not an eigenvalue of $(-\Delta)_p^s$.

Assume by contradiction that $\lambda^*=\infty$. According to $(F_1)$, there exists $\widehat{\lambda}>1$ such that $(P_{\widehat{\lambda}})$ has a bounded solution $u_{\widehat{\lambda}}$, and $\widehat{\lambda}f(t)>\mu t^{p-1}$ for all $t>0$, so $u_{\widehat{\lambda}}$ is a supersolution of the following equation
\begin{equation}\label{2209141646}
(-\Delta)_p^s u = \mu |u|^{p-2}u \;\text{ in}~\Omega,\quad u =0~\;\; \text{in}~ \mathbb{R}^N\setminus\Omega.
\end{equation}
Furthermore, let $\varphi_1>0$ be the eigenfunction associated to $\lambda_1$ (which is known to be continuous). We choose $\varepsilon > 0$ small such that
$$\lambda_1(\varepsilon\varphi_1(x))^{p-1}\le \widehat{\lambda} f(0)\le \widehat{\lambda} f(u_{\widehat{\lambda}}(x)) \quad \mbox{for all $x\in \Omega$.}$$ By comparison principle, we get $\varepsilon\varphi_1\le u_{\widehat{\lambda}}$. Using once more Lemma \ref{Super_Sub_Le},
%
there exists a nontrivial solution to \eqref{2209141646}, which is absurd seeing the choice of $\mu$. This contradiction shows that $0<\lambda^*<\infty$. 

Next, for $\lambda\in (0, \lambda^*)$, we claim that $u_\lambda$ is a stable solution. It follows from Proposition \ref{2308061322} and $u_{\lambda}\in L^{\infty}(\Omega)$ that there exists a stable solution $0<\widehat{u}_{\lambda}\le u_\lambda$. Since $u_{\lambda}$ is the minimal solution, we get $\widehat{u}_{\lambda}={u}_{\lambda}$. So ${u}_{\lambda}$ is a stable solution.

For any $0<\lambda'<\lambda''<\lambda^*$, $u_{\lambda''}$ is a supersolution of $(P_{\lambda'})$, which yields that $u_\lambda$ is increasing with respect to $\lambda$.
\end{proof}

\subsection{Dichotomy result for $W_0^{s, p}(\Omega)$ solutions}
Here $f$ satisfies $(F_1)$ and $f^{\frac{1}{p-1}}$ is convex, we will show that $\lambda^*$ is also the dichotomy parameter for existence of $W_0^{s, p}(\Omega)$ solution. We mainly follow the idea in \cite{BC1996}.

Let $\varepsilon\in (0, 1)$, and 
\[
g=f^{\frac{1}{p-1}},\quad \widetilde{g}_\e =(1-\varepsilon)^{\frac{1}{p-1}}g; \quad 
h(u)=\int_0^{u}\frac{1}{g(t)}dt,\quad \widetilde{h}_\e(u)=\int_0^{u}\frac{1}{\widetilde{g}_\e(t)}dt=\frac{h(u)}{(1-\varepsilon)^{\frac{1}{p-1}}}.
\]
Denote $\Phi_\e(u)=\widetilde{h}_\e^{-1}\circ h$. By convexity of $g$ and \cite[Lemma 4]{BC1996}, $\Phi_\e\in C^2(\mathbb{R}_+)$ and we have 
\begin{itemize}
\item[--] $\Phi_\e(0)=0$ and $0\le \Phi_\e(u)\le u$ for all $u\ge 0$, $\e \in (0, 1)$;
\item[--] $\Phi_\e$ is increasing, concave and $\Phi_\e'(u)\le 1$ in $\R_+$;
\item[--] If $h(\infty)< \infty$, then $\Phi_\e(\infty)<\infty$.
\end{itemize}

\begin{proof}[\bf Proof of Theorem \ref{2312030111}]  Assume that $(P_\lambda)$ has a solution $u_\lambda\in W_0^{s,p}(\Omega)$. For any $\varepsilon\in (0, 1)$, put $v_{\lambda}:=\Phi_\e(u_\lambda)$. Applying Proposition \ref{2308041528} with $\Psi=-\Phi_\e$, and test function space as $W_0^{s,p}(\Omega)\cap L^{\infty}(\Omega)$, we have
\begin{equation}\label{2312192051}
(-\Delta)_p^sv_{\lambda}=(-\Delta)_p^s \Phi_\e(u_\lambda)\ge \lambda \Phi_\e'(u_\lambda)^{p-1}f(u_\lambda)=\lambda (1-\varepsilon) f(v_{\lambda}).
\end{equation}
Note that $W_0^{s,p}(\Omega)\cap L^{\infty}(\Omega)$ is dense in $W_0^{s,p}(\Omega)$, and $f(v_{\lambda})\in W_0^{s,p}(\Omega)^*$ since $0 \le f(v_\lambda)\le f(u_\lambda)$. Thus $v_{\lambda}$ is a supersolution of $(P_{(1-\varepsilon)\lambda})$. We will consider two subcases:

{\it Case $h(\infty)<\infty$}. $v_{\lambda} \in L^\infty(\Omega)$. By Lemma \ref{Super_Sub_Le}, there is a bounded solution for $(P_{(1-\varepsilon)\lambda})$.

{\it Case $h(\infty)=\infty$}. By the concavity of $h$, there holds
\begin{equation}\label{2312192044}
h(u_\lambda)\le h(v_{\lambda})+\frac{u_\lambda-v_{\lambda}}{f(v_{\lambda})^{\frac{1}{p-1}}}.
\end{equation}
Combining $h(v_{\lambda})=(1-\varepsilon)^{\frac{1}{p-1}}h(u_\lambda)$ with $\varepsilon \in (0, 1)$, we deduce that
\begin{equation}\label{2312022206}
\Big(1-(1-\varepsilon)^{\frac{1}{p-1}}\Big) f(v_{\lambda})^{\frac{1}{p-1}}\le \frac{u_\lambda}{h(u_\lambda)}\le C(1+u_\lambda).
\end{equation}
Here we use an iterative procedure similar to that for Lemma \ref{2306271914}, we only consider briefly the case $N>sp$, and omit the detail for other situations. The Sobolev embedding $W_0^{s,p}(\Omega)\subset L^{q_0}(\Omega)$ and \eqref{2312022206} imply $f(v_{\lambda})\in L^{\frac{q_0}{p-1}}(\Omega)$, where $q_0=p_s^*$ as $N>sp$. Then Proposition \ref{2306271234} leads to
$$v_{\lambda}\in L^{q_1}(\Omega),\quad q_1=\frac{Nq_0}{N-\frac{spq_0}{p-1}}.$$
Recall that $v_{\lambda}$ is a supersolution of $(P_{(1-\varepsilon)\lambda})$, so there exists a weak solution $w_1$ of $(P_{(1-\varepsilon)\lambda})$ such that $0\le w_1\le v_{\lambda}$. By the same procedure as above, we get
\[
w_1\in L^{q_2}(\Omega),\;\; q_2=\frac{Nq_1}{N-\frac{spq_1}{p-1}}, \quad \mbox{if }\; q_1 < \frac{N(p-1)}{sp}.
\]
Repeating the iteration, whenever $q_n < \frac{N(p-1)}{sp}$, we get $q_{n+1}=\frac{Nq_n}{N-\frac{spq_n}{p-1}} > q_n$, and $w_n \in L^{q_{n+1}}(\Omega)$ solution to $(P_{(1-\varepsilon)^{n}\lambda})$. As before, there is $n_0 \geq 1$ (independent on $\e$) such that $q_{n_0} >\frac{N(p-1)}{sp}$, therefore $(P_{(1-\varepsilon)^{n_0}\lambda})$ has a solution $w_{n_0}\in L^{\infty}(\Omega)$, which means that $(1-\varepsilon)^{n_0}\lambda \le \lambda^*$. As $\varepsilon$ is arbitrary, we obtain $\lambda \le \lambda^*$. \end{proof}

\section{Regularity of extremal function for power-like nonlinearities}\label{2312021519}
We consider now the regularity of the extremal function $u^*= \lim_{\lambda\to \lambda^*}u_\lambda$ via uniform estimates for stable solutions to \eqref{2306182230}. Here we handle nonlinearities with power-like growth, i.e. $g$ satisfies \eqref{2306202233} and \eqref{2306211610}. 

Applying Lemma \ref{2312192311} with $h(t)=|t|^{2\gamma-2}t$ $(\gamma>1)$, $p=2$ and $H(t) = \frac{\sqrt{2\gamma-1}}{\gamma}|t|^{\gamma-1}t$, we get
\begin{equation}\label{2306202030}
\left(a^{\gamma}-b^{\gamma}\right)^2\le \frac{\gamma^2}{2\gamma-1}(a-b)(a^{2\gamma-1}-b^{2\gamma-1}), \quad \forall \; a, b\ge 0, \gamma>1.
\end{equation}

\begin{proof}[\bf Proof of Theorem \ref{2306211611}]
Assume that $u$ is a non-negative stable solution of \eqref{2306182230}. Let  $T_k(t) = \min\{k, |t|\}$ ($k > 0$),
and $\psi=uT_k(u)^{\gamma-1}$ for $\gamma>1$. Clearly $\psi\in \mathcal{T}_u(\Omega)$. By the stability assumption \eqref{2306201531},
\begin{equation}\label{2306202053}
\int_{\Omega}g_t(x, u)\psi^2 dx\le (p-1)\int_{\mathbb{R}^{2N}}{\frac{w(x, y)}{|x-y|^{N+ps}}}dxdy
\end{equation}
where
\[
w(x, y)=|u(x)-u(y)|^{p-2}(uT_k(u)^{\gamma-1}(x)-uT_k(u)^{\gamma-1}(y))^2.
\]
We divide $\mathbb{R}^{2N}$ into four parts:
\[
\begin{aligned}
&A_1:=\{u(x)>k,~ u(y)>k\},\quad A_2:=\{u(x)\le k,~ u(y)\le k\},\\
&A_3:=\{u(x)\le k,~ u(y)>k\},\quad A_4:=\{u(x)>k,~ u(y)\le k\}.
\end{aligned}
\]
Direct computations yield
\[
\begin{aligned}
I_1:=&\int_{A_1}{\frac{w(x, y)}{|x-y|^{N+ps}}}dxdy
= k^{2\gamma-2}\int_{A_1}{\frac{|u(x)-u(y)|^{p}}{|x-y|^{N+ps}}}dxdy,\\
I_2:=&\int_{A_2}{\frac{w(x, y)}{|x-y|^{N+ps}}}dxdy= \int_{A_2}{\frac{|u(x)-u(y)|^{p-2}(u^{\gamma}(x)-u^{\gamma}(y))^2}{|x-y|^{N+ps}}}dxdy,\\
I_3:=&\int_{A_3\cup A_4}{\frac{w(x, y)}{|x-y|^{N+ps}}}dxdy=2\int_{A_3}{\frac{|u(x)-u(y)|^{p-2}(u^{\gamma}(x)-k^{\gamma-1}u(y))^2}{|x-y|^{N+ps}}}dxdy.
\end{aligned}
\]
Using $\varphi=uT_k(u)^{2\gamma-2}$ as test function to \eqref{2306182230}, we have
\begin{equation}\label{2306202054}
\int_{\Omega}g(x, u)uT_k(u)^{2\gamma-2}dx=\int_{\mathbb{R}^{2N}}{\frac{h(x, y)}{|x-y|^{N+ps}}}dxdy,
\end{equation}
where
\[
h(x, y)=|u(x)-u(y)|^{p-2}(u(x)-u(y))(uT_k(u)^{2\gamma-2}(x)-uT_k(u)^{2\gamma-2}(y)).
\]
Therefore
\[
\begin{aligned}
H_1:=&\int_{A_1}{\frac{h(x, y)}{|x-y|^{N+ps}}}dxdy
= k^{2\gamma-2}\int_{A_1}{\frac{|u(x)-u(y)|^{p}}{|x-y|^{N+ps}}}dxdy,\\
H_2:=&\int_{A_2}{\frac{h(x, y)}{|x-y|^{N+ps}}}dxdy= \int_{A_2}{\frac{|u(x)-u(y)|^{p-2}(u(x)-u(y))(u^{2\gamma-1}(x)-u^{2\gamma-1}(y))}{|x-y|^{N+ps}}}dxdy,\\
H_3:=&\int_{A_3\cup A_4}{\frac{h(x, y)}{|x-y|^{N+ps}}}dxdy
=2\int_{A_3}{\frac{|u(x)-u(y)|^{p-2}(u(x)-u(y))(u^{2\gamma-1}(x)-k^{2\gamma-2}u(y))}{|x-y|^{N+ps}}}dxdy.
\end{aligned}
\]
Obviously $I_1 = H_1$. By the inequality \eqref{2306202030}, for $\gamma > 1$, 
\begin{equation}\label{2306202049}
I_2\le \frac{\gamma^2}{2\gamma-1}H_2.
\end{equation}
Moreover, for $(x, y)\in A_3$, applying again \eqref{2306202030},
\[
\begin{aligned}
 \frac{2\gamma-1}{\gamma^2}(u^{\gamma}(x)-k^{\gamma-1}u(y))^2& \le \left[u(x)-k^{\frac{\gamma-1}{\gamma}}u(y)^{\frac1\gamma}\right] \times \Big[u(x)^{2\gamma-1}-k^{\frac{(2\gamma-1)(\gamma-1)}{\gamma}}u(y)^{\frac{2\gamma-1}\gamma}\Big]\\
& = \big[u(x)-u(y)\big]\times \Big[u(x)^{2\gamma-1}-k^{2\gamma-2}u(y)\Big]\\
&\quad -\Big[u(x)u(y)-k^{\frac{\gamma-1}{\gamma}}u(x)u(y)^{\frac{1}{\gamma}}\Big]\times \Big[k^{\frac{(2\gamma-1)(\gamma-1)}{\gamma}}u(y)^{\frac{\gamma-1}{\gamma}}-u(x)^{2\gamma-2})\Big]\\
& \le  \big[u(x)-u(y)\big]\times \Big[u(x)^{2\gamma-1}-k^{2\gamma-2}u(y)\Big].
\end{aligned}
\]
This means that for $\gamma > 1$,
\begin{equation}\label{2306202050}
I_3\le \frac{\gamma^2}{2\gamma-1}H_3.
\end{equation}
Remark that $\frac{\gamma^2}{2\gamma-1} > 1$ for any $\gamma > 1$. By \eqref{2306202053}-\eqref{2306202050},  there holds
\begin{equation}\label{2306202051}
\int_{\Omega}g_t(x, u)u^2T_k(u)^{2\gamma-2} dx\le \frac{(p-1)\gamma^2}{2\gamma-1}\int_{\Omega}g(x, u)uT_k(u)^{2\gamma-2} dx.
\end{equation}
On the other hand, it follows from \eqref{2306202233}, \eqref{2306211610} that for any $\overline{m}\in (p-1, m)$, there is $A\ge 0$ such that
\[
g_t(x, t)t\ge \overline{m}g(x,t)-\overline{m}C(1+A)^m\chi_{\{t\le A\}} \mbox{ for all }t\ge 0 \mbox{ and a.e. }x\in \Omega,
\]
where $\chi$ is the characteristic function. Combining with \eqref{2306202051}, we get
\begin{equation}\label{2310241701}
\left(\overline{m}-\frac{(p-1)\gamma^2}{2\gamma-1}\right)\int_{\Omega}g(x, u)uT_k(u)^{2\gamma-2}dx\le \overline{m}C(1+A)^mA^{2\gamma-1}|\Omega|.
\end{equation}
Take 
\begin{equation}\label{2310241859}
\gamma\in \left(\frac{\overline{m}}{p-1}, \beta(\overline m)\right) \quad \mbox{with }\; \beta(\sigma) := \frac{\sigma+\sqrt{\sigma(\sigma-p+1)}}{p-1} \;\mbox{  for }\; \sigma \geq p-1, 
\end{equation}
so that $\overline{m}>\frac{(p-1)\gamma^2}{2\gamma-1}$.
Using \eqref{2310241701} and Fatou's lemma, let $k \to \infty$, we arrive at
\begin{equation}\label{2306202234}
\int_{\Omega}g(x, u)u^{2\gamma-1}dx\le C,
\end{equation}
where $C$ depends only on $p$, $\Omega$, $A$, $\gamma$ and $g$. Coming back to the equation \eqref{2306182230}, there holds
\[
\begin{aligned}
\|u\|_{W_0^{s, p}(\Omega)}^p = \int_{\Omega}g(x, u)u dx & \le \int_{\{u\le 1\}}g(x, u) dx+\int_{\{u>1\}}g(x, u)u^{2\gamma-1} dx \le C'
\end{aligned}
\]
which yields the uniform boundedness of stable solutions to \eqref{2306182230} in $W_0^{s, p}(\Omega)$.

Furthermore, combining \eqref{2306202234} with \eqref{2306202233} and letting $\overline m \to m$, we see that
\begin{equation}\label{2306202238}
\int_{\Omega}|g(x, u)|^qdx\le C \quad\text{for}~1\le q<\frac{2\beta(m)+m-1}{m}.
\end{equation}
To establish $u \in L^\infty(\Omega)$, from \eqref{2306202238} and Lemma \ref{2306271914}, it is enough to have
\[
\left(1-\frac{p-1}{m}\right)\frac{N}{sp}<\frac{2\beta(m)+m-1}{m},
\]
or equivalently, $N<G_s(m, p)$. The proof is done.
\end{proof}


As explained, Theorem \ref{2307200048} is somehow direct consequence of Theorem \ref{2306211611}.
\begin{proof}[\bf Proof of Theorem \ref{2307200048}]
As $m>p-1$ and $f$ satisfies \eqref{newThm1.9}, we can apply Theorem \ref{2306211611} with $g(x, t) = \lambda f(t)$ and $\lambda < \lambda^*$. It follows that $\{u_\lambda\}$ is bounded in $W_0^{s,p}(\Omega)$. Hence
\begin{equation}\label{2307182203}
u_\lambda\rightharpoonup u^* \text{~in~}W_0^{s,p}(\Omega).
\end{equation}
On the other hand, we claim $f( u^*)\in W_0^{s,p}(\Omega)^*$. For example,  when $N > sp$, we have $\|\lambda f(u_\lambda)\|_{(p_s^*)'} dx<C$. Indeed, applying \eqref{2306202238} with $g(x, u) = \lambda f(u)$, $\lambda < \lambda^*$, it suffices to show $$\frac{Np}{Np-N+sp}<\frac{2\beta({m})+{m}-1}{{m}},$$
which is valid by
\[
\frac{Np}{Np-N+sp}<\frac{p}{p-1}<\frac{\beta{(m)}+m}{m}<\frac{2\beta({m})+{m}-1}{{m}}.
\]
By $f(u^*)\in W_0^{s,p}(\Omega)^*$ and monotone convergence, passing to the limit $\lambda\to \lambda^*$, there holds
\[
\langle (-\Delta)_p^s u^*, \varphi \rangle=\lambda^*\int_{\Omega}f(u^*)\varphi dx, \quad\forall\, \varphi\in W_0^{s,p}(\Omega), \; \varphi \geq 0.
\]
Then working with positive and negative part for general $W_0^{s,p}(\Omega)$ function, we see that $u^*$ is a solution of $(P_{\lambda^*})$. We can claim that $u^*$ is a minimal solution of $(P_{\lambda^*})$, and $u^*$ is stable as in the proof of Theorem \ref{2306211344}. So Claim (i) is proved. Claim (ii) can be proven similarly, we omit the detail.
\end{proof}

\section{Regularity of extremal solution for convex nonlinearities}\label{2312021527}
Here we consider $(P_\lambda)$ with the assumptions $(F_1)$-$(F_2)$, that permit to work with more general nonlinearities such as $f(t)=e^t$. For convenience, we denote 
\begin{align}
\label{Def6.0}
\widetilde{f}(t)=f(t)-f(0), \;\; \psi(t)=\widetilde{f}(t)^{\frac{1}{p-1}}, \;\; \varphi(t)=\int_0^t\psi'(\sigma)^2d\sigma\;\; \mbox{and}\;\; \zeta(t)=\int_0^t\psi'(\sigma)(\psi'(t)-\psi'(\sigma))d\sigma.
\end{align}
We follow some ideas from \cite{Nedev2000, Sanchon_PA2007, Sanchon07}. Since the approach for $N\le sp$ is very similar to the case $N>sp$, without loss of generality, we will {\bf assume} $N>sp$ in this section.

\subsection{Proof of Theorem  \ref{2309111457}}
Let $\lambda \in (0, \lambda^*)$, as $u_\lambda \in W_0^{s, p}(\Omega)\cap L^\infty(\Omega)$, $\psi(u_\lambda) \in {\mathcal T}_{u_\lambda}(\Omega)$. Take $\psi(u_\lambda)$ as test function in the stability condition \eqref{2306201531}, we have
\begin{equation}\label{2309111531}
\lambda \int_{\Omega}f'(u_\lambda) \psi(u_\lambda)^2 dx\le (p-1)\int_{\mathbb{R}^{2N}}{\frac{|u(x)-u(y)|^{p-2}(\psi(u_\lambda(x))-\psi(u_{\lambda}(y)))^2}{|x-y|^{N+ps}}}dxdy.
\end{equation}
Take $\varphi(u_\lambda)$ as test function to $(P_\lambda)$, there holds
\begin{equation}\label{2309111556}
\int_{\mathbb{R}^{2N}}{\frac{J_{u_\lambda}(x, y)(\varphi(u_\lambda(y))-\varphi(u_\lambda(x)))}{|x-y|^{N+ps}}}dxdy=\lambda\int_{\Omega}f(u_\lambda)\varphi(u_\lambda)dx.
\end{equation}
From \eqref{Def6.0}-\eqref{2309111556} and Cauchy-Schwarz's inequality,
\begin{align*}
(\psi(x)-\psi(y))^2&=\Big(\int_{u_\lambda(x)}^{u_\lambda(y)}\psi'(t)dt\Big)^2\le \big[u_\lambda(y)-u_\lambda(x)\big]\times\big[\varphi(u_\lambda(y))-\varphi(u_\lambda(x))\big],
\end{align*}
therefore
\begin{equation}\label{2310242209}
 \int_{\Omega}\widetilde{f}'(u_\lambda) \psi(u_\lambda)^2 dx\le (p-1)\left(\int_{\Omega}\widetilde{f}(u_\lambda)\varphi(u_\lambda)dx+f(0)\int_{\Omega}\varphi(u_\lambda)dx\right).
\end{equation}
It is easy to check that
$\widetilde{f}'(t)\psi(t)^2-(p-1)\widetilde{f}(t)\varphi(t)=(p-1)\widetilde{f}(t)\zeta(t)$, 
from \eqref{2310242209} we get 
\begin{equation}\label{2309132002}
\int_{\Omega}\psi(u_\lambda)^{p-1}\zeta(u_\lambda) dx\le f(0)\int_{\Omega}\varphi(u_\lambda) dx.
\end{equation}

On the other hand, as $\psi$ is convex and superlinear at infinity, there holds
\begin{equation}\label{2309131710}
\underset{t\to\infty}{\lim}\frac{\zeta(t)}{\psi'(t)}=\infty.
\end{equation}
By the convexity of $\psi$, we have
\begin{align}
\label{new6.1}
\int_{\Omega}\frac{\psi(u_\lambda)^p}{u_\lambda} dx \le \int_{\Omega}\psi(u_\lambda)^{p-1}\psi'(u_\lambda) dx.
\end{align}
Using \eqref{2309132002} and \eqref{2309131710}, there exists $C > 0$ such that
\begin{align}
\label{new6.2}
\begin{split}
\int_{\Omega}\psi(u_\lambda)^{p-1}\psi'(u_\lambda) dx \leq C + \frac{1}{2f(0)}\int_{\Omega}\psi(u_\lambda)^{p-1}\zeta(u_\lambda) dx & \le C + \frac{1}{2}\int_{\Omega}\varphi(u_\lambda) dx\\
& \le C' +  \frac{1}{2}\int_{\Omega}\psi(u_\lambda)^{p-1}\psi'(u_\lambda) dx.
\end{split}
\end{align}
For the last inequality, we used $p\ge 2$ and 
\[
\varphi(t)\le \psi'(t)\int_0^t\psi'(\sigma)d\sigma = \psi(t)\psi'(t)\le \psi(t)^{p-1}\psi'(t) \quad \mbox{for $t$ large such that $\psi(t) \geq 1$}.
\]
From \eqref{new6.1} and \eqref{new6.2}, there holds 
\begin{align}
\label{new6.3}
\int_\Omega \frac{\widetilde{f}(u_\lambda)^{p'}}{u_\lambda}dx = \int_{\Omega}\frac{\psi(u_\lambda)^p}{u_\lambda} dx \le \int_{\Omega}\psi(u_\lambda)^{p-1}\psi'(u_\lambda) dx \leq C.
\end{align} 
Applying \eqref{new6.3}, we can claim: 
\begin{align}
\label{newc1}
\mbox{If $\{u_\lambda\}$ is bounded in $L^q(\Omega)$ for $q \geq 1$, then $\{f(u_\lambda)\}$ is bounded in $L^{\frac{qp'}{q+1}}(\Omega)$.}
\end{align}
In fact, denote 
\[
\Omega_1:=\{x\in \Omega: \widetilde{f}(u_\lambda(x))^{\frac{p'}{q+1}}\ge u_\lambda(x) \}\;\;
\text{and}\;\;
\Omega_2:=\{x\in \Omega: \widetilde{f}(u_\lambda(x))^{\frac{p'}{q+1}}< u_\lambda(x) \}.
\]
By \eqref{new6.3}, there holds
\begin{align*}
\int_{\Omega}\widetilde{f}(u_\lambda)^{\frac{qp'}{q+1}} dx&=\int_{\Omega_1}\widetilde{f}(u_\lambda)^{\frac{qp'}{q+1}} dx+\int_{\Omega_2}\widetilde{f}(u_\lambda)^{\frac{qp'}{q+1}} dx \le \int_{\Omega_1}\frac{\widetilde{f}(u_\lambda)^{p'}}{u_\lambda}dx+ \int_{\Omega_2}u_{\lambda}^{q}dx\le C.
\end{align*}

By $(F_1)$ and \eqref{new6.3}, we deduce that $\{f(u_\lambda)\}$ is bounded in $L^1(\Omega)$. 
Proposition \ref{2306271234} implies that $\{u_\lambda\}$ is bounded in $L^q(\Omega)$ for all $1\le q<q_0:=\frac{N(p-1)}{N-sp}$. By induction with the claim \eqref{newc1} and Proposition \ref{2306271234}, we get an increasing sequence $\{q_k\}$ such that $\{u_\lambda\}$ is bounded in $L^{q_k}(\Omega)$ and $\{f(u_\lambda)\}$ is bounded in $L^\frac{q_kp'}{q_k+1}(\Omega)$, where 
\begin{equation}
\label{newite1}
q_{k+1}=\frac{Nq_k(p-1)p'}{N(q_k+1)-spq_kp'},\quad \mbox{ whenever }\; N >\frac{spq_kp'}{q_k+1}.
\end{equation}
Therefore, by simple analysis, we will meet two cases. 
\begin{itemize}
\item[--] If {$N<spp'$}, there is $k_0 \geq 1$ such that $N > \frac{spq_{k_0}p'}{q_{k_0}+1}$. Therefore $\{u_\lambda\}$ is bounded in $L^{\infty}(\Omega)$;
\item[--] If {$N\ge spp'$}, then $\lim_{k\to \infty}q_k = \frac{N(p-1)}{N - spp'}$, $\{u_\lambda\}$ is bounded in $L^{q}(\Omega)$ for all $1\le q < \frac{N(p-1)}{N-spp'}$.
\end{itemize}
The corresponding integral estimates for $f(u_\lambda)$ is given by \eqref{newc1}. 
Furthermore, we remark that $f(u_\lambda)u_\lambda$ is bounded in $L^1(\Omega)$, hence $\{u_\lambda\}$ is bounded in $W_0^{s,p}(\Omega)$. We conclude the corresponding regularity results for $u^*$ by monotone and weak convergence argument as before, so we omit the details.

\subsection{Proof of Theorem \ref{2307201716}}
Since $\psi$ is convex and $\psi(0)=0$, we have $\psi'(t)\ge\frac{\psi(t)}{t}$ and $\psi'(t)\ge \psi'(\sigma)$  for any $t\ge \sigma>0$, equivalently
\begin{equation}\label{2407142339}
f'(t)\ge (p-1)\frac{\widetilde{f}(t)}{t}\;\; \mbox{and} \;\; f'(\sigma)\le f'(t)\left(\frac{\widetilde{f}(t)}{\widetilde{f}(\sigma)}\right)^{\frac{2-p}{p-1}}\quad\forall\; t\ge \sigma>0.
\end{equation}
Next, for any $\gamma\ge \frac{1}{p-1}$, we define the following function
\[
\kappa(t):=\int_0^t \widetilde{f}(\sigma)^{2\gamma-2}f'(\sigma)^2d\sigma.
\]

\begin{proposition}\label{2308071331}
Assume $(F_1)$-$(F_2)$. If there exists $\gamma\ge \frac{1}{p-1}$ such that 
\begin{equation}\label{2407142322}
\underset{t\to\infty}{\lim\sup} \frac{(p-1)\gamma^2\kappa(t)}{\widetilde{f}(t)^{2\gamma-1}f'(t)}<1,
\end{equation}
then $u^*\in W_0^{s,p}(\Omega)$ is a stable solution of $(P_{\lambda^*})$. Moreover:
\begin{itemize}
\item[\rm (i)] If $N<(2\gamma+1)sp$, then $u^*\in L^{\infty}(\Omega)$;
\item[\rm (ii)] If $N\ge (2\gamma+1)sp$, then $u^*\in L^{q}(\Omega)$ for all $q\in [1, q_*)$, where
\begin{equation}\label{2312011717}
q_*=
\begin{cases}
\begin{aligned}
&\frac{N[(p-1)(2\gamma+1)-1]}{N-(2\gamma+1)sp} \;\; &&\mbox{if}\; N>(2\gamma+1)sp;\\
&\infty \;\; &&\mbox{if}\; N=(2\gamma+1)sp.
\end{aligned}
\end{cases}
\end{equation}
\end{itemize}
\end{proposition}
\begin{proof}
Let $\eta:=\widetilde{f}(u_\lambda)^{\gamma}$, then $\eta \in\mathcal{T}_{u_\lambda}(\Omega)$ since $u_\lambda\in W_0^{s, p}(\Omega)\cap L^{\infty}(\Omega)$ and $\gamma\ge \frac{1}{p-1}$. Using H\"older inequality, it holds that
\[
\begin{aligned}
(\eta(x)-\eta(y))^2=\gamma^2\left(\int_{u_\lambda(x)}^{u_\lambda(y)}\widetilde{f}(t)^{\gamma-1}f'(t)dt\right)^2\le\gamma^2 [u_\lambda(y)-u_\lambda(x)]\times [\kappa(u_\lambda(y))-\kappa(u_\lambda(x))].
\end{aligned}
\]
Combining the above inequality and the stability condition \eqref{2306201531} with $\eta$, we get
\begin{equation}\label{2307131020}
\lambda \int_{\Omega}\widetilde{f}(u_\lambda)^{2\gamma}f'(u_\lambda) dx\le (p-1)\gamma^2\int_{\mathbb{R}^{2N}}{\frac{J_{u_\lambda}(x, y)(\kappa(u_\lambda(y))-\kappa(u_\lambda(x)))}{|x-y|^{N+ps}}}dxdy.
\end{equation}
Multiplying $(P_\lambda)$ by $\kappa(u_\lambda)$, we see that
\begin{equation}\label{2307131021}
\int_{\mathbb{R}^{2N}}{\frac{J_{u_\lambda}(x, y)(\kappa(u_\lambda(y))-\kappa(u_\lambda(x)))}{|x-y|^{N+ps}}}dxdy=\lambda\int_{\Omega}f(u_\lambda)\kappa(u_\lambda)dx,
\end{equation}
which combines with \eqref{2307131020} deduces
\begin{equation}\label{2307131123}
\begin{aligned}
\int_{\Omega}\widetilde{f}(u_\lambda)^{2\gamma}f'(u_\lambda) dx&\le (p-1)\gamma^2\left(\int_{\Omega}\widetilde{f}(u_\lambda)\kappa(u_\lambda)dx+f(0)\int_{\Omega}\kappa(u_\lambda)dx\right).
\end{aligned}
\end{equation}
By \eqref{2407142322} and $(F_1)$, there hold 
\begin{equation}\label{2307131125}
\underset{t\to \infty}{\lim \sup}(p-1)\gamma^2\frac{\widetilde{f}(t)\kappa(t)}{\widetilde{f}(t)^{2\gamma}f'(t)}<1\quad\mbox{and} \quad \underset{t\to \infty}{\lim}(p-1)\gamma^2\frac{\kappa(t)}{\widetilde{f}(t)^{2\gamma}f'(t)}=0.
\end{equation}
By \eqref{2407142339}, \eqref{2307131123}-\eqref{2307131125},  we can derive
\begin{equation}\label{2307131710}
\int_{\Omega}\frac{\widetilde{f}(u_\lambda)^{2\gamma+1}}{u_\lambda} dx \leq C \int_{\Omega}\widetilde{f}(u_\lambda)^{2\gamma}f'(u_\lambda) dx\le C'.
\end{equation}
Proceed exactly as in the proof of Theorem \ref{2309111457}, we see that 
\begin{align}
\label{newc2}
\mbox{If $\{u_\lambda\}$ is bounded in $L^q(\Omega)$ for $q \geq 1$, then $\{f(u_\lambda)\}$ is bounded in $L^{\frac{q(2\gamma +1)}{q+1}}(\Omega)$.}
\end{align}
As before,  there is a sequence $\{q_k\}$ such that $\{u_\lambda\}$ is bounded in $L^{q_k}(\Omega)$ and $\{f(u_\lambda)\}$ is bounded in $L^\frac{q_k(2\gamma+1)}{q_k+1}(\Omega)$, where
\begin{equation*}
q_{k+1}=\frac{Nq_k(p-1)(2\gamma +1)}{N(q_k+1)-spq_k(2\gamma +1)},\quad \mbox{ if }\; N >\frac{spq_k(2\gamma+1)}{q_k+1}.
\end{equation*}
Whenever the bootstrap argument can start, that is $q_0 < q_1$, we will meet two cases. 
\begin{itemize}
\item[--] If {$N<sp(2\gamma +1)$}, $\{u_\lambda\}$ and $\{f(u_\lambda)\}$ are bounded in $L^{\infty}(\Omega)$;
\item[--] If {$N\ge sp(2\gamma +1)$}, $\{u_\lambda\}$ is bounded in $L^{q}(\Omega)$ for all $1\le q <q_*$.
\end{itemize}

Suppose $N\ge (2\gamma +1)sp$. In virtue of \eqref{newc2}, $\{f(u_\lambda)\}$ is bounded in $L^{q}(\Omega)$ for all $1\le q< \ell$ with
\[
\ell =\frac{q_*(2\gamma +1)}{q_*+1} = \frac{N((2\gamma+1)(p-1)-1)}{N(p-1)-sp}>\frac{Np}{N(p-1)+sp}.
\]
The last inequality holds since $\gamma\ge\frac{1}{p-1}$. Hence $f(u_\lambda)\to f(u^*)$ in $W_0^{s,p}(\Omega)^*$. In all these cases, we can conclude that $u^*$ is a stable solution of $(P_{\lambda^*})$.
\end{proof}

The Proposition above tells us that when $\gamma\ge\frac{1}{p-1}$ such that \eqref{2407142322}, we will get a higher regularity for $u^*$. In fact, the condition \eqref{2308071234} for $\tau$  guarantees the existence of such $\gamma$. The following two lemmas (see \cite[Lemmas 12-13]{Sanchon07}) give more precisely the corresponding range of $\gamma$ and the growth of $f$ for given $\tau$.
\begin{lemma}
\label{2308071334}
Let $\tau$ be given in  \eqref{2308071234}. If $\tau>\frac{p-2}{p-1}$, then \eqref{2407142322} holds for any
\[
\gamma\in \left(\frac1{p-1}, \frac{1+\sqrt{1-(p-1)(1-\tau)}}{p-1} \right).
\]
\end{lemma}
\begin{lemma}\label{2308072041}
Assume that $f\in C^2(\mathbb{R}^+, \mathbb{R}^+)$ such that $f'(t)>0$ for $t>0$. Let $\tau$ be given in  \eqref{2308071234}. If $\tau<1$, then for any $\varepsilon\in (0, 1-\tau)$, there exists $c=c(\varepsilon)>0$ such that
\[
f(t)\le c(1+t)^{\frac{1}{1-(\tau+\varepsilon)}}\;\;\mbox{for all }t\ge 0.
\]
\end{lemma}
\begin{proof}[\bf Proof of Theorem \ref{2307201716} completed] Since $\tau>\frac{p-2}{p-1}$, using Proposition \ref{2308071331} and Lemma \ref{2308071334}, we have that $u^*\in W_0^{s, p}(\Omega)$ is a stable solution of $(P_{\lambda^*})$. Moreover, if
\[
N<N(s, p):=sp+\frac{2sp}{p-1}\left(1+\sqrt{1-(p-1)(1-\tau)}\right),
\]
then $u^*\in L^{\infty}(\Omega)$, and if $N\ge N(s, p)$, then $u^*\in L^q(\Omega)$ for all $1\le q<q_0$, where
\[
q_0=\frac{N\left(2\sqrt{1-(p-1)(1-\tau)}+p\right)}{N-N(s, p)}.
\]

{\sl Case $\tau\ge 1$}. By direct computations, $N(s, p)>sp+\frac{4sp}{p-1}$, so when $N<sp+\frac{4sp}{p-1}$, we  have $u^*\in L^{\infty}(\Omega)$.

{\sl Case $\frac{p-2}{p-1}<\tau< 1$}. For any $\varepsilon\in (0, 1-\tau)$, applying Lemma \ref{2306271914} and Lemma \ref{2308072041} with $m=\frac{1}{1-(\tau+\varepsilon)}$, for the sake of $u^*\in L^{\infty}(\Omega)$, it suffices to have
\[
q_0>\frac{N}{sp}\left(1-\frac{p-1}{m}\right),
\]
which is equivalent to
\[
N<\frac{sp}{p-1}\left(1+\frac{p}{1-(p-1)(1-\tau-\varepsilon)}+\frac{2\sqrt{1-(p-1)(1-\tau)}}{1-(p-1)(1-\tau-\varepsilon)}\right).
\]
By $\tau>\frac{p-2}{p-1}$, one can check that
\[
\frac{sp}{p-1}\left(1+\frac{p}{1-(p-1)(1-\tau-\varepsilon)}+\frac{2\sqrt{1-(p-1)(1-\tau)}}{1-(p-1)(1-\tau-\varepsilon)}\right)>sp+\frac{4sp}{p-1}.
\]
Therefore, whenever $N<sp+\frac{4sp}{p-1}$, one has $u^*\in L^{\infty}(\Omega)$.
\end{proof}

\noindent \textbf{Data availability} Data sharing is not applicable to this article as no new data
were created or analyzed in this study.

\medskip
\noindent \textbf{Acknowledgement} The author completed this paper during his Ph.D. in  East China Normal University. The author would like to thank his supervisor  Dong Ye for careful reading of the manuscript.

\end{document}